\documentclass[12pt,reqno,a4paper]{amsart}
\usepackage{amsmath}
\usepackage{amssymb}
\usepackage{stmaryrd}
\usepackage{amsthm}

\usepackage[applemac]{inputenc}

\advance\voffset    by -1  cm
\advance\hoffset    by -1.5cm
\advance\textwidth  by  3  cm
\advance\textheight by  1  cm
\sffamily

\usepackage{pstricks}
\usepackage{multido}

\usepackage{graphicx}
\usepackage{color}
\newcommand{\showgrid}{}

\newcommand{\gridon}{\renewcommand{\showgrid}{\psset{subgriddiv=1,griddots=10,gridlabels=6pt}\psgrid}}

\gridon % enable grid if uncommented

\newgray{hellgrau}{0.9}
\newrgbcolor{dummycolor}{0.75 0.85 0.0}

\long\def\comment#1{%
\psshadowbox[%
fillstyle=solid,fillcolor=dummycolor,linewidth=0.05%
]{\hbox to 12.70cm {\vbox{\strut {\bf Kommentar:} #1\hfill}}}%
}

\def\figref#1{Figure~\ref{#1}}
\def\heute{\number\day.~\ifcase\month\or
 J\"anner\or Februar\or M\"arz\or April\or Mai\or Juni\or
 Juli\or August\or September\or Oktober\or November\or Dezember\fi
 \space\number\year}

\def\defeq{:=}
\def\absof#1{\left|#1\right|}
\def\setof#1{\left\{#1\right\}}
\def\of#1{\!\left(#1\right)}

\def\pas#1{\left(#1\right)}
\def\brk#1{\left[#1\right]}
\def\sgn{{\operatorname{sgn}}}

\def\floor#1{{\left\lfloor #1\right\rfloor}}

\def\bit{\begin{itemize}}
\def\eit{\end{itemize}}
\def\beq{\begin{equation}}
\def\eeq{\end{equation}}

\def\e{{\mathbf e}}
\def\1{{\mathbf 1}}

{\begin{list}{#1}{\parsep0em \itemsep0.3em \labelwidth1em
\labelsep0.5em \leftmargin1.5em }} {\end{list}}

\def\path{{%\mathrm 
P}}
\def\figref#1{Figure~\ref{#1}}

\def\EM#1{{\em #1\/}}

\newif\ifenglish
\englishtrue

\newtheorem{thm}{\ifenglish Theorem\else Satz\fi}
\newtheorem{pro}{Proposition}
\newtheorem{lem}{Lemma}
\newtheorem{rem}{\ifenglish Remark\else Bemerkung\fi}
\newtheorem{dfn}{Definition}
\newtheorem{cor}{\ifenglish Corollary\else Korollar\fi}

\theoremstyle{remark}
\newtheorem{obs}{\ifenglish Observation\else Beobachtung\fi}

\def\secref#1{Section~\ref{#1}}

\def\pfaffian{\operatorname{Pf}}

\def\inv{\operatorname{inv}}
\def\matchings#1{{\mathcal M}_{#1}}
\def\superpositions#1{{\mathcal S}_{#1}}
\def\nofmatchings#1{M_{#1}}
\def\subgraph#1#2{\brk{#1-#2}}
\def\subgraphedges#1#2#3{\brk{#1-#2-#3}}
\def\weight{\omega}
\def\vertices{{\mathrm V}}
\def\edges{{\mathrm E}}
\def\orientation{\xi}
\def\symbolf{f\!}

\def\concat{{\scriptstyle\circ}}
\def\word#1{\pas{#1}}
\def\complement#1{\overline{#1}}

\def\symdef#1#2{#1\,{\vartriangle}\,#2} % symmetric difference
\def\numof#1{\#\of{#1}}                % number of a set
\def\cardof#1{\left\vert{#1}\right\vert}                % cardinality of a set
\def\smatch{\mu}                        % single matching
\def\bicoloured#1#2#3{#1_{{ #2}\vert{ #3}}}
\def\theblues{{\blue\mathbf{b}}}
\def\thereds{{\red\mathbf{r}}}
\def\thecoloureds{{\black\mathbf{c}}}
\def\thewhites{{\black\mathbf{w}}}
\def\swapcolours#1{\chi_{#1}}
\def\SETSUM#1#2{{{\mathbf\Sigma}\of{{#1\subseteq#2}}}}
\def\symm{\mathfrak{S}}
\def\descof#1{d\,\of{#1}}

\def\RHS{right--hand side}
\def\symmeven#1{{\symm^{0}_{#1}}}
\def\step#1{Step #1: }
\def\semibip#1#2{S_{#1,#2}}
\def\completeG#1{K_{#1}}

\def\word#1{{\mathbf #1}}
\def\pmatch{matching}
\def\Pmatch{Matching}
\def\strangeset{{\mathcal F}\of{G\|X}}
\def\markextra#1{#1}
\def\disjuni{\,\dot{\cup}\,}

\def\and{\wedge}
\def\indexed#1#2{{{\mathbf #1}}_{#2}}

\def\complement#1{\overline{#1}}
\def\thelefts{{\red{\mathbf R}}}
\def\therights{{\blue{\mathbf B}}}
\def\myrange#1{\setof{1,\dots,k}}
\def\som{superposition of matchings}
\def\soms{superpositions of matchings}
\def\nofcooriented#1{{\vert #1\vert}^{\circlearrowright}}
\def\nofinterior#1{{\vert #1\vert}^{\circ}}

\newcount\labelcount

\newif\iflongversion
\longversiontrue

\newif\ifverylongversion
\verylongversionfalse

\def\loopless{\relax}
\def\wpit{involution--type}

\begin{document}

\bibliographystyle{plain}

\title[Graphical condensation: Overlapping Pfaffians]{Graphical condensation, overlapping Pfaffians and superpositions of matchings}

\begin{abstract}
The purpose of this note is to exhibit clearly how
the ``graphical condensa\-tion'' identities of Kuo, Yan, Yeh and Zhang follow from
classical Pfaffian identities by the Kasteleyn--Percus method for the
enumeration of \pmatch s. Knuth termed the relevant identities
``overlapping Pfaffian'' identities and the key concept of proof ``superpositions
of matchings''. In our uniform presentation of the material, we also give an apparently
unpublished general ``overlapping Pfaffian'' identity of Krattenthaler.

A previous version of this paper contained an erroneous application of
the Kasteleyn--Percus method, which is now corrected. 
\end{abstract}

\author{Markus Fulmek}
\address{Fakult\"at f\"ur Mathematik, % \\ Universit\"at Wien\\
Nordbergstra\ss e 15, A-1090 Wien, Austria}
\email{{\tt Markus.Fulmek@Univie.Ac.At}\newline\leavevmode\indent
{\it WWW}: {\tt http://www.mat.univie.ac.at/\~{}mfulmek}
}

\date{\today}
\thanks{
Research supported by the National Research Network ``Analytic
Combinatorics and Probabilistic Number Theory'', funded by the
Austrian Science Foundation. 
}

\maketitle

%\nocite{fulmek:2001,kasteleyn:1967,knuth:1996,kuo:2006,percus:1969,yan-zhang:2005,mb:minc:1978}

%$\pfaffian$ is similar to $\det$ for skew--symmetric matrices
%of even size, with permutations replaced by perfect matchings.
%$\per$ is ``$\det$ without the signs''. $\hafnian$ is similar
%to $\pfaffian$ for symmetric matrices, without signs.

% \input todo

\section{Introduction}

%\bit
%\item Consider finite ordered set of vertices; order is inherited in subsets,
%complete graph on set thus implicitly directed?: For our purposes, any simple edge--weighted
%graph on this set can also be viewed as complete graph, where some edge weights are zero.
%\item Allow ``double edges''.
%\item orientation of graphs.
%\eit

%\alert{Explain: Why present known stuff once more?}
%\bit
%\item Lange Geschichte: Ideen aus dem 19.~Jahrhundert,
%\item Vielfach wiederentdeckt,
%\item Neufassung kombinatorischer Behandlung von ``Overlapping Pfaffians''
%	im Sinne von Knuth und Hamel,
%\item Kleine Neuigkeit (Krattenthalers Identität)
%\eit

%Planarer Graph: In der Ebene zeichnen - graphical! Daher auch viele Illustrationen.

%Achtung: Matchings = Perfect matching!

%Vollständiger graph -- Notation $K_n$, $K_V$.

%Alle Mengen endlich.

In the last 5 years, several authors
\cite{kuo:2006,kuo:2004,propp:2003,yan-zhang:2005,yan-yeh-zhang:2005} % seemingly rediscovered (special cases of)
came up with identities related to the enumeration of matchings in planar graphs, together
with a beautiful method of proof, which they termed \EM{graphical condensation}.

In this paper, we show that these identities are special cases of certain \EM{Pfaffian identities} (in the simplest case Tanner's identity \cite{tanner:1878}), by simply
applying the Kasteleyn--Percus method \cite{kasteleyn:1963,percus:1969}. These identities involve products of Pfaffians, for which Knuth
\cite{knuth:1996} coined the term \EM{overlapping Pfaffians}. Overlapping Pfaffians were further investigated by Hamel \cite{hamel:2001}.

Knuth gave a very clear and concise exposition not only of the
results, but also of the main idea of proof, which he termed \EM{superposition of
matchings}.

Tanner's identity dates back to the 19th century --- and so does the basic idea of superposition of matchings, which was used for a proof of Cayley's Theorem \cite{cayley:1849}
by Veltmann in 1871 \cite{veltmann:1871} and independently by Mertens in 1877
\cite{mertens:1877} (as was already pointed out by Knuth \cite{knuth:1996}). Basically
the same proof of Cayley's Theorem was presented by Stembridge \cite{stembridge:1990},
who gave a very elegant ``graphical'' description of Pfaffians.

The purpose of this note is to exhibit clearly how ``graphical condensation'' is % tightly connected to
 connected to ``overlapping Pfaffian'' identities. This is achieved by
\bit
\item using Stembridge's description of Pfaffians to give a simple, uniform
  presentation of the underlying idea of
 ``superposition of matchings'', accompanied by many graphical illustrations (which
 should demonstrate \EM{ad oculos} the beauty of this idea),
 \item using this idea to give \EM{uniform proofs} for several known ``overlapping Pfaffian'' identities
 	and 
% \item presenting
 a general ``overlapping Pfaffian'' identity,
 which to the best of our knowledge  is due to
 Krattenthaler \cite{krattenthaler:2006} and was not published before,
\item and (last but not least) making clear how the ``graphical condensation'' identities
 of Kuo \cite[Theorem~2.1 and Theorem~2.3]{kuo:2006}, Yan, Yeh and Zhang
 \cite[Theorem~2.2 and Theorem~3.2]{yan-yeh-zhang:2005} and Yan and Zhang \cite[Theorem~2.2]{yan-zhang:2005} follow immediately from the
 ``overlapping Pfaffian'' identities  via the classical Kasteleyn--Percus method
 for the enumeration of (perfect) matchings.
\eit

This note is organized as follows:
\bit
\item \secref{sec:graphs-matchings} presents the basic definitions and notations
	used in this note,
\item \secref{sec:superpositions} presents the concept of ``superposition of matchings'',
	using a simple instance of ``graphical condensation'' as an introductory example,
%	and explains how ``graphical condensation'' makes use of it,
\item \secref{sec:pfaffians} presents Stembridge's description of
Pfaffians\iflongversion\ together with the ``superposition of matchings''--proof of Cayley's Theorem\fi,
\item \secref{sec:kasteleyn-percus} recalls %(quite briefly)
	the Kasteleyn--Percus method,
\item \secref{sec:tanner} presents Tanner's classical identity and more general
 ``overlapping Pfaffian'' identities, together with ``superposition of matchings''--proofs,
 and deduces the ``graphical condensation'' identities
 \cite[Theorem~2.1 and Theorem~2.3]{kuo:2006}, 
 \cite[Theorem~2.2 and Theorem~3.2]{yan-yeh-zhang:2005} and \cite[Theorem~2.2]{yan-zhang:2005}.
 % for them.
\eit

% \subsection{Notational conventions}

\section{Basic notation: Ordered sets (words), graphs and \Pmatch s}
\label{sec:graphs-matchings}
The \EM{sets} we shall consider in this paper will always be
% \bit 
% \item
\EM{finite}
% \item
and \EM{ordered},
% \eit
whence we may view them as \EM{words} of distinct letters
$$
\word{\alpha} = \setof{\alpha_1,\alpha_2,\dots,\alpha_n} \simeq
\pas{\alpha_1,\alpha_2,\dots,\alpha_n}.
$$
%If $\word\beta = \pas{\beta_1,\beta_2,\dots,\beta_m}$ is another set (word) with
%$\word\alpha\cap\word\beta=\emptyset$, then by $\word{\alpha\beta}$ we mean the
%\EM{concatenation of words}
%$$
%\word{\alpha\beta}\defeq\pas{\alpha_1,\dots,\alpha_n,\beta_1,\dots,\beta_n}.
%$$
When considering a subset $\word\gamma\subseteq\word\alpha$, we shall \EM{always} assume that
the elements (letters) of $\word\gamma$ appear in the same order as in $\word\alpha$, i.e.,
$$
\gamma=\setof{\alpha_{i_1},\alpha_{i_2},\dots,\alpha_{i_k}}
\simeq\pas{\alpha_{i_1},\alpha_{i_2},\dots,\alpha_{i_k}} \text{ with } i_1<i_2<\dots<i_k.
$$

We choose this somewhat indecisive notation because the order of the elements (letters)
is not always relevant. For instance, for \EM{graphs} $G$ we shall employ the
usual (set--theoretic) terminology:  $G=G\of{\vertices,\edges}$ with vertex set
$\vertices\of G = \vertices$ and edge set $\edges\of G=\edges$, and the ordering
of $V$ is irrelevant for typical graph--theoretic questions like ``is $G$ a \EM{planar}
graph?''.

The graphs we shall
consider in this paper will always be \EM{finite} and \EM{loopless} (they may, however,
have multiple edges). Moreover, the graphs will always be \EM{weighted}, i.e., we assume
a weight function $\weight:\edges\of G\to R$, where $R$
is some integral domain.
(If we are interested in mere enumeration, we may simply choose $\weight\equiv 1$.)

The \EM{weight} $\weight\of U$ of some \EM{subset} of edges
$U\subseteq \edges\of G$ is defined as
$$\weight\of U\defeq\prod_{e\in U}\weight\of e.$$

The \EM{total weight} (or \EM{generating function}) of some
\EM{family} $\mathcal F$ of subsets of $\edges\of G$ is defined as
$$\weight\of{\mathcal F}\defeq\sum_{U\in\mathcal F}\weight\of U.$$

For some subset $S\subseteq \vertices\of G$, we denote by $\subgraph{G}{S}$ the
subgraph of $G$ \EM{induced} by the vertex set $\pas{\vertices\of G\setminus S}$.

A \EM{\pmatch} in $G$ is a subset $\smatch\subseteq \edges\of G$
of edges such that
\bit
\item
no two edges in $\smatch$ have a vertex in common,
\item
and every vertex in $\vertices\of G$ is incident with precisely one
edge in $\smatch$.
\eit
(This concept often is called a \EM{perfect \pmatch}).
Note that a \pmatch\ $\smatch$ may be viewed as a \EM{partition} of
$\vertices\of{G}$ into blocks of cardinality 2.

\section{Kuo's Proposition and superposition of \pmatch s}
\label{sec:superpositions}

Denote the family of \EM{all} \pmatch s of $G$ by $\matchings{G}$, and
denote the total weight of all \pmatch s of $G$ by
$\nofmatchings{G}\defeq \weight\of{\matchings{G}}$.

According to Kuo \cite{kuo:2006}, the following proposition is
a generalization of results of Propp \cite[section 6]{propp:2003} and Kuo
\cite{kuo:2004}, and was first proved combinatorially by
Yan, Yeh and Zhang \cite{yan-yeh-zhang:2005}:
\begin{pro}
\label{prop:kuo}
Let $G$ be a \EM{planar} graph with four vertices $a$, $b$, $c$ and $d$ that
appear in that cyclic order on the boundary of a face of $G$. Then
\begin{equation}
\label{eq:kuos-proposition}
\nofmatchings{G}\nofmatchings{\subgraph{G}{\setof{a,b,c,d}}} +
\nofmatchings{\subgraph{G}{\setof{a,c}}}\nofmatchings{\subgraph{G}{\setof{b,d}}} =
% \\
\nofmatchings{\subgraph{G}{\setof{a,b}}}\nofmatchings{\subgraph{G}{\setof{c,d}}} +
\nofmatchings{\subgraph{G}{\setof{a,d}}}\nofmatchings{\subgraph{G}{\setof{b,c}}}.
\end{equation}
\end{pro}

As we will see, this statement is a direct consequence of Tanner's \cite{tanner:1878} identity
(see \cite[Equation (1.0)]{knuth:1996}) and the Kasteleyn--Percus method
\cite{kasteleyn:1967}, but we shall use it here as a simple example to
introduce the concept of \EM{superposition of matchings}, as the
% We present a simple proof for this statement which arises from a
straightforward combinatorial intepretation of the products
involved in equations like \eqref{eq:kuos-proposition} was termed by Knuth \cite{knuth:1996}.
%coined the term \EM{superposition of matchings}.
% Following Knuth \cite{knuth:1996}, we call
%this  \EM{superpositions of matchings} (the name was given
%by Knuth \cite{knuth:1996}):
% of total weights: % appearing in \eqref{eq:kuos-proposition}.

\subsection{Superpositions of matchings}
\begin{figure}
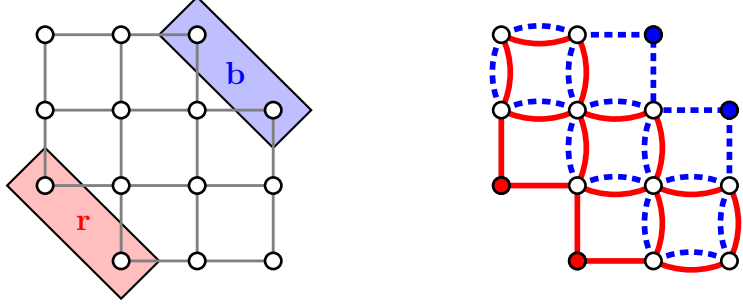

\caption{Illustration: A graph $G$ with two disjoint subsets of vertices $\thereds$
and $\theblues$ (shown in the left picture) gives rise to the bicoloured graph
$\bicoloured{G}{\thereds}{\theblues}$ (shown in the right picture; blue edges are shown
as dashed lines).}
\label{fig:bicoloured}
\input graphics/bicoloured
\end{figure}

%\begin{dfn}[bicoloured graph]
Consider a simple graph $G$ and two disjoint (but otherwise arbitrary) subsets of vertices
$\theblues\subseteq\vertices\of G$ 
and $\thereds\subseteq\vertices\of G$. Call $\theblues$ the \EM{blue} vertices,
$\thereds$ the \EM{red} vertices, $\thecoloureds\defeq\thereds\cup\theblues$ the 
\EM{coloured} vertices and the remaining $\thewhites\defeq\vertices\of G\setminus\thecoloureds$ the \EM{white} vertices.
Now consider the  \EM{bicoloured} graph
$B=\bicoloured{G}{\thereds}{\theblues}$
\bit
\item with vertex set $\vertices\of B\defeq\vertices\of G$,
	% \setminus \pas{\theblues\cap \thereds}$,
\item and with edge set $\edges\of B$ equal to the \EM{disjoint union} of 
	\bit
	\item the edges of $\edges\of{\subgraph{G}{\theblues}}$, which are coloured \EM{red},
	\item and the edges of $\edges\of{\subgraph{G}{\thereds}}$, which are coloured \EM{blue}.
	\eit
\eit
Here, ``disjoint union'' should be understood in the sense
that $\edges\of{\subgraph{G}{\theblues}}$ and $\edges\of{\subgraph{G}{\thereds}}$
are subsets of two different ``copies'' of $\edges\of{G}$, respectively.
This concept will appear frequently in the following: assume that we have two copies
of some set $M$. We may imagine these copies to have different colours, red and blue,
and denote them accordingly by $M_r$ and $M_b$, respectively. Then ``by definition'' subsets
$A\subseteq M_r$ and $B\subseteq M_b$ are \EM{disjoint}: every element in $A\cap B$ (in the
ordinary sense, as subsets of $M$) appears twice (as red copy and as blue copy)
in $A\cup B$. Introducing the notation $X\disjuni Y$
as a shortcut for \EM{disjoint union}, i.e, for ``$X\cup Y$, where $X\cap Y=\emptyset$'', we can write:
$$\edges\of B=\edges\of{\subgraph{G}{\theblues}}\disjuni\edges\of{\subgraph{G}{\thereds}}.$$
%then any two sets $A$ and $B$
%In general, for two subsets $A,B\subseteq M$ of some set $M$
%we introduce the notation
%$$A\disjuni B$$
%meaning that $A$ and $B$ are to be interpreted as subsets of different copies of $M$,
%whence (in this peculiar sense) $A\cap B=\emptyset$.
%We may imagine these copies of $M$ to have different colours, red and blue, whence elements
%in  $A\cap B$ (in the ordinary sense) appear twice (as red copy and as blue copy)
%in $A\disjuni B$.
Note that in $B=\bicoloured{G}{\thereds}{\theblues}$
\bit
\item all edges incident with blue vertices (i.e., with vertices in $\theblues$)
	are blue,
\item all edges incident with red vertices (i.e., with vertices in $\thereds$)
	are red,
\item and all edges in $\edges\of{\subgraph{G}{\thecoloureds}}$ appear as
\EM{double} edges in
$\edges\of B$; one coloured red and the other coloured blue.
\eit
See \figref{fig:bicoloured}
for an illustration.

%where all the edges in $\edges\of{\subgraph{G}{U}}$ are coloured red,
%and all the edges in $\edges\of{\subgraph{G}{W}}$ are coloured blue;
The weight function $\weight$ on the edges of graph $B=\bicoloured{G}{\thereds}{\theblues}$ 
is assumed to be inherited from graph $G$: $\weight\of e$ in $B$ equals $\weight\of e$ in
$G$ (irrespective of the \EM{colour} of $e$ in $B$).
%\end{dfn}

\begin{obs}[superposition of matchings]
\label{obs:interpretation}
Define the weight $\weight\of{X,Y}$ of a \EM{pair} $\pas{X,Y}$ of subsets of edges as
$$\weight\of{X,Y}\defeq\weight\of X\cdot\weight\of Y.$$%=\weight\of{X\disjuni Y}.$$ 
Then
$\nofmatchings{\subgraph{G}{\theblues}}\nofmatchings{\subgraph{G}{\thereds}}$
clearly equals the total weight of
$\matchings{\subgraph{G}{\theblues}}\times\matchings{\subgraph{G}{\thereds}}$; and
the typical summand in % the total weight
$\nofmatchings{\subgraph{G}{\theblues}}\nofmatchings{\subgraph{G}{\thereds}}$
is the product of the weights $\weight\of{{\red\smatch_r}}\cdot\weight\of{{\blue\smatch_b}}$
of some pair of \pmatch s
$\pas{{\red\smatch_r},{\blue\smatch_b}}\in\matchings{\subgraph{G}{\theblues}}\times\matchings{\subgraph{G}{\thereds}}$. Such pair of \pmatch s can be viewed as the 
disjoint union ${\red\smatch_r}\disjuni{\blue\smatch_b}\subseteq \edges\of B$ in the bicoloured graph $B$, where ${\red\smatch_r}$ is a subset of the \EM{red} edges, and ${\blue\smatch_b}$ is a subset of
the \EM{blue} edges. We call any subset in $\edges\of B$
%From this point of view,
%$\nofmatchings{\subgraph{G}{U}}\nofmatchings{\subgraph{G}{W}}$
%appears as the total weight of all bicoloured subsets of $\edges\of B$
which arises from a pair of \pmatch s in this way a
\EM{superposition of matchings}, and we
denote by $\superpositions{B}$ the family of superpositions
of matchings of $B$. So there is a weight preserving bijection
\begin{equation}
\label{eq:translation}
\matchings{\subgraph{G}{\theblues}}\times\matchings{\subgraph{G}{\thereds}}
\leftrightarrow
\superpositions{B}.
\end{equation}
\end{obs}

\begin{obs}[nonintersecting bicoloured paths/cycles]
\label{obs:nonintersecting}
It is obvious that some subset $S\subseteq \edges\of B$
of edges of the bicoloured graph $B$ is a superposition of matchings if and only if
\bit
\item every blue vertex $v$ (i.e., $v\in\theblues$)
   %(call these vertices the \EM{blue vertices})
	is incident with precisely one
	blue edge from $S$,
\item every red vertex $v$ (i.e., $v\in\thereds$) %(call these vertices the \EM{red vertices})
	is incident with precisely one
	red edge from $S$,
\item every white vertex $v$ (i.e., $v\in\thewhites$)
   %in $\vertices\of B \setminus \pas{{\theblues}\cup{\thereds}}$
	% (call these vertices the \EM{white vertices})
	is incident
	with precisely one blue edge \EM{and} precisely one red edge from $S$.
\eit
Stated otherwise: A superposition of matchings in $B$ may be viewed as a % \EM{nonintersecting}
family of
% \EM{bicoloured}
paths and cycles,
\bit
\item such that every vertex of $B$ belongs to \EM{precisely one}
	path or cycle (i.e., the paths/cycles are \EM{nonintersecting}: no two
	different cycles/paths have a vertex in common),
\item such that edges of each cycle/path \EM{alternate} in colour along the cycle/path
	(therefore, we call them \EM{bicoloured}: Note that a bicoloured cycle must have
	even length),
\item such that precisely the \EM{end vertices} of paths are \EM{coloured} (i.e., red or
	blue), and all other vertices are white.
\eit
\end{obs}

Note that a bicoloured cycle of length $>2$ in the bicoloured graph
$B=\bicoloured{G}{\thereds}{\theblues}$ corresponds to a cycle in $G$,
while a bicoloured cycle of length $2$ in $B$ corresponds to a ``doubled edge''
in $G$.

\begin{obs}[colour--swap along paths]
\label{obs:swap}
For an arbitrary \EM{coloured} vertex $x$ in some superposition
of matchings $S$ of $\edges\of{B}$, we may \EM{swap} colours for all the edges in the
unique path $p$ in $S$ with end vertex $x$ (see \figref{fig:recolouring}). Without loss of
generality, assume that ${\red x}$ is red. Depending on the colour of
the other end vertex $y$ of $p$, this colour--swap results
in a set of coloured edges $\overline S$, which is a
superposition of matchings in
\bit
\item $B^\prime=\bicoloured{G}{{\thereds^\prime}}{{\theblues^\prime}}$, where
	$\thereds^\prime\defeq\pas{\thereds\setminus\setof{{\red x}}}\cup\setof{{\blue y}}$
	and
	$\theblues^\prime\defeq\pas{\theblues\setminus\setof{{\blue y}}}\cup\setof{{\red x}}$,
	if ${\blue y}$ is blue (i.e., of the opposite colour as ${\red x}$;
	the length of the path $p$ is \EM{even} in this case --- this case is illustrated
	in \figref{fig:recolouring}),
\item $B^{\prime\prime}=\bicoloured{G}{{\thereds^{\prime\prime}}}{{\theblues^{\prime\prime}}}$, where
	$\thereds^{\prime\prime}\defeq\pas{\thereds\setminus\setof{{\red x, y}}}$
	and
	$\theblues^{\prime\prime}\defeq\theblues\cup\setof{{\red x, y}}$,
	if ${\red y}$ is red (i.e., of the same colour as ${\red x}$;
	the length of path $p$ is \EM{odd} in this case).
\eit
Clearly, this operation of swapping colours defines a \EM{weight preserving}
injection 
\begin{equation}
\swapcolours{{x}}:
\superpositions{B}
\to
\pas{
\superpositions{B^\prime}
\disjuni
\superpositions{B^{\prime\prime}}
}
\end{equation}
(which, viewed as mapping onto its image, is an involution: $\swapcolours{{x}}=\swapcolours{{x}}^{-1}$).
So $\swapcolours{{x}}$ together with the bijection \eqref{eq:translation} gives
a weight preserving injection
\begin{equation}
\matchings{\subgraph{G}{\theblues}}\times\matchings{\subgraph{G}{\thereds}}
\to
\bigcup_{\pas{B^\prime,B^{\prime\prime}}}
	\matchings{\subgraph{G}{\theblues^\prime}}
	\times\matchings{\subgraph{G}{\thereds^\prime}}
\disjuni
\matchings{\subgraph{G}{\theblues^{\prime\prime}}}
	\times\matchings{\subgraph{G}{\thereds^{\prime\prime}}},
\end{equation}
where the union is over all pairs $\pas{B^\prime,B^{\prime\prime}}$
of bicoloured graphs that arise by the recolouring of the path $p$, as described
above.
\end{obs}

\begin{figure}
\caption{Illustration: Take graph $G$ of \figref{fig:bicoloured} and consider a
matching in $\subgraph{G}{\thereds}$ ($\thereds=\red\setof{x,t}$), whose edges
are colored \EM{blue} (shown as dashed lines), and a matching in $\subgraph{G}{\theblues}$ ($\theblues=\blue\setof{y,z}$), whose edges
are colored \EM{red}. This superposition of matchings determines
a unique path $p$ connecting ${\red x}$ and ${\blue y}$ in the bicoloured graph $\bicoloured{G}{\thereds}{\theblues}$. Swapping the colours of the edges of $p$ determines uniquely a matching in
$\subgraph{G}{\thereds^\prime}$ ($\thereds^\prime=\red\setof{y,t}$) and a matching
in $\subgraph{G}{\theblues^\prime}$ ($\theblues^\prime=\blue\setof{x,z}$).}
\label{fig:recolouring}
\input graphics/recolouring
\end{figure}

\subsection{The ``graphical condensation method''}
\label{sec:graphical-condensation}
Now we apply the reasoning outlined in Observations~\ref{obs:interpretation},
\ref{obs:nonintersecting} and \ref{obs:swap} for the proof of
Proposition~\ref{prop:kuo} (basically the same proof is presented in \cite{kuo:2006}):
\begin{proof}[Proof of Proposition~\ref{prop:kuo}]
Clearly, for all the superpositions of \pmatch s
(see Observation~\ref{obs:interpretation})
% four products of total weights
involved in % (the combinatorial interpretation of)
\eqref{eq:kuos-proposition}, the set
$\thecoloureds$ of coloured vertices in the associated bicoloured graphs is
$\setof{a,b,c,d}$. In any superposition of \pmatch s,
there are two \EM{nonintersecting} paths (see Observation~\ref{obs:nonintersecting})
with end vertices
in $\setof{a,b,c,d}$. Since $G$ is \EM{planar} and the vertices
$a$, $b$, $c$ and $d$ appear in this cyclic order in the boundary of a \EM{face} $F$ of $G$, 
the path starting in vertex $a$ cannot end in vertex $c$ (otherwise
it would intersect the path connecting $b$ and $d$; see \figref{fig:start-graphics}
for an illustration).

So consider the
bicoloured graphs
\bit
\item
$B_1\defeq\bicoloured{G}{\thereds_1}{\theblues_1}$ with $\thereds_1\defeq\setof{a,b,c,d}$, $\theblues_1=\emptyset$,
\item and
$B_2\defeq\bicoloured{G}{\thereds_2}{\theblues_2}$ with $\thereds_2\defeq\setof{a,c}$, $\theblues_2=\setof{b,d}$.
\eit
%and set
%\bit
%\item $B_1\defeq\bicoloured{G}{\thereds_1}{\theblues_1}$
%\item and $B_2\defeq\bicoloured{G}{\thereds_2}{\theblues_2}$.
%\eit
Observe that
\begin{multline*}
\nofmatchings{G}\nofmatchings{\subgraph{G}{\thereds_1}} +
\nofmatchings{\subgraph{G}{\theblues_2}}\nofmatchings{\subgraph{G}{\thereds_2}} = \\
\weight\pas{
	\pas{\matchings{\subgraph{G}{\theblues_1}}\times\matchings{\subgraph{G}{\thereds_1}}}
	\disjuni
	\pas{\matchings{\subgraph{G}{\theblues_2}}\times\matchings{\subgraph{G}{\thereds_2}}}
} =
%	\\
\weight\of{
	\superpositions{B_1}\disjuni
	\superpositions{B_2}
}.
\end{multline*}

%bicoloured graphs
%%\bit
%%\item
%$B_1\defeq B\of{\emptyset,\setof{a,b,c,d}}$
%% , corresponding to
%%	$\matchings{G}\times\matchings{\subgraph{G}{a,b,c,d}}$
%% \item
%and $B_2\defeq B\of{\setof{a,c},\setof{b,d}}$:
%% , corresponding to
%%	$\matchings{\subgraph{G}{a,c}}\matchings{\subgraph{G}{b,d}}$
%% \eit
%%and observe that the total weight of all superposition subsets
%%in these two graphs equals the left hand
%%side of \eqref{eq:kuos-proposition}:
%% Observation~\ref{obs:interpretation} gives a weight--preserving
%% bijection
%\begin{multline}
%\nofmatchings{G}\nofmatchings{\subgraph{G}{a,b,c,d}} +
%\nofmatchings{\subgraph{G}{a,c}}\nofmatchings{\subgraph{G}{b,d}} = \\
%\weight\pas{
%	\matchings{G}\times\matchings{\subgraph{G}{a,b,c,d}} \cup
%	\matchings{\subgraph{G}{a,c}}\times\matchings{\subgraph{G}{b,d}}
%} =
%%	\\
%\weight\of{
%	\superpositions{B\of{\emptyset,{a,b,c,d}}}\cup
%	\superpositions{B\of{\setof{a,c},\setof{b,d}}}
%}
%\end{multline}
Note that for any superposition of \pmatch s, the other end--vertex of the
bicoloured path starting at $a$ necessarily has
\bit
\item the \EM{same colour} as $a$ in $B_1$ (i.e., \EM{red}),
\item the \EM{other colour} as $a$ in $B_2$ (i.e., \EM{blue}).
\eit
(See \figref{fig:matchings}.) So consider the
bicoloured graphs
\bit
\item
$B_1^\prime\defeq\bicoloured{G}{\thereds_1^\prime}{\theblues_1^\prime}$ with $\thereds_1^\prime\defeq\setof{b,c}$, $\theblues_1^\prime=\setof{a,d}$,
\item and
$B_2^\prime\defeq\bicoloured{G}{\thereds_2^\prime}{\theblues_2^\prime}$ with $\thereds_2^\prime\defeq\setof{c,d}$, $\theblues_2^\prime=\setof{a,b}$.
\eit
It is easy to see that the operation $\swapcolours{a}$ of
swapping colours of edges along the path starting at vertex $a$
(see Observation~\ref{obs:swap})
defines a weight preserving involution
\begin{equation*}
	\swapcolours{a}:
	\superpositions{B_1}\disjuni
	\superpositions{B_2}
\leftrightarrow
	\superpositions{B_1^\prime}\disjuni
	\superpositions{B_2^\prime},
\end{equation*}
and thus gives a weight preserving involution
\begin{multline*}
	\pas{\matchings{G}\times\matchings{\subgraph{G}{\thereds_1}}}
		\disjuni
	\pas{
		\matchings{\subgraph{G}{\theblues_2}}\times
		\matchings{\subgraph{G}{\thereds_2}}
	}
\leftrightarrow \\
	\pas{
		\matchings{\subgraph{G}{\theblues_1^\prime}}\times
		\matchings{\subgraph{G}{\thereds_1^\prime}}
	}
		\disjuni
	\pas{
		\matchings{\subgraph{G}{\theblues_2^\prime}}\times
		\matchings{\subgraph{G}{\thereds_2^\prime}}
	}.
\end{multline*}
(See \figref{fig:matchings} for an illustration.)
%
%Now, for every superposition subset
%in $B_1$ as well as in $B_2$, swap edge--colours along the path
%starting in $a$.  It is easy to see, that the resulting 
%
%\bit
%\item $\nofmatchings{\subgraph{G}{a,b}}\nofmatchings{\subgraph{G}{c,d}}$
%\item $\nofmatchings{\subgraph{G}{a,d}}\nofmatchings{\subgraph{G}{b,c}}$
%\eit
%
%For an arbitrary subset $U\subseteq \edges\of G$, a perfect matching of $\subgraph{G}{U}$ constitutes a matching of $G$ .
%Therefore, for arbitrary subsets $U$, $W$ of $\edges\of G$, we may interpret the product
%$\nofmatchings{\subgraph{G}{U}}\nofmatchings{\subgraph{G}{W}}$
%as the total weight of \EM{superpositions of matchings} of
%$G$ in the following sense. For every pair of
%perfect matchings
%$\pas{\bar{U}, \bar{W}}$, $\bar{U}\subseteq{\edges\of{\subgraph{G}{U}}}$
%and $\bar{W}\subseteq{\edges\of{\subgraph{G}{W}}}$, the union
%$\bar{U}\cup\bar{W}$ is a multiset of edges whose weight
%$\weight\of{\bar{U}\cup\bar{W}}$ is in one--to--one--correspondence to a
%summand of $\nofmatchings{\subgraph{G}{U}}\nofmatchings{\subgraph{G}{W}}$.
%%, which is the total weight of all such ``superpositions of
%%matchings''.
%
\end{proof}

\begin{figure}
\caption{A simple planar graph $G$ with vertices $a$, $b$, $c$ and $d$ appearing in this
order in the boundary of face $F$.}
\label{fig:start-graphics}
\input graphics/start-graphics
\end{figure}

\begin{figure}
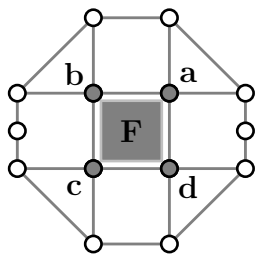

\caption{Take the planar graph $G$ from \figref{fig:start-graphics} and consider the
bicoloured graphs $B_1$, $B_2$, $B_1^\prime$ and $B_2^\prime$ from the proof of Proposition~\ref{prop:kuo}: The pictures show
certain superpositions of \pmatch s in these bicoloured graphs (the edges belonging
to the \pmatch s are drawn as thick lines, the blue edges appear as
dashed lines). The arrows indicate the possible effects of the operation $\swapcolours{a}$,
i.e., of swapping colours of edges in the unique path with end vertex $a$.}
\label{fig:matchings}
\input graphics/matchings
\end{figure}

This bijective proof % for Proposition~\ref{prop:kuo} 
certainly is very
satisfactory. But since there is a
well--known powerful method for enumerating perfect matchings
in planar graphs, namely the Kasteleyn--Percus method (see
\cite{kasteleyn:1963,kasteleyn:1967,percus:1969}) which involves \EM{Pfaffians},
the question arises whether Proposition~\ref{prop:kuo} (or the bijective method
of proof) gives additional insight or provides a new perspective.

%For reader's convenience, we shall repeat the corresponding
%concepts and facts.

\section{Pfaffians}
\label{sec:pfaffians}
The name \EM{Pfaffian} was introduced by Cayley \cite{cayley:1852} (see
\cite[page 10f]{knuth:1996} for a concise historical survey).
Here, we follow closely Stembridge's exposition \cite{stembridge:1990}:
\begin{dfn}
\label{dfn:Stembridge}
Consider the \EM{complete graph} $\completeG{V}$ on the (ordered) set of vertices
$V=\pas{v_1, \dots v_n}$, with weight function
$\weight:\edges\of{\completeG{V}}\to R$.
Draw this graph in the upper halfplane in the following specific way:
\bit
\item Vertex $v_i$ is represented by the point $\pas{i,0}$,
\item edge $\setof{v_i,v_j}$ is represented by the half--circle with center
	$\pas{\frac{i+j}{2},0}$ and radius $\frac{\absof{j-i}}{2}$.
\eit
(See the left picture in \figref{fig:Kn-matchings}).

Consider some \pmatch\ $\mu=\setof{\setof{v_{i_1},v_{j_1}},\dots,\setof{v_{i_m},v_{j_m}}}$ in $\completeG V$.
Clearly, if such $\mu$ exists, then $n=2 m$ must be even. By convention, we assume
${i_k}<{j_k}$ for $k=1,\dots,m$. A \EM{crossing} of $\mu$ corresponds to a crossing
of edges in the specific drawing just described, or more formally: A crossing of $\mu$ is a
pair of edges $\pas{\setof{v_{i_k},v_{j_k}}, \setof{v_{i_l},v_{j_l}}}$ of $\mu$ such that
$$
i_k < i_l < j_k < j_l.
$$
Then the sign of $\mu$ is defined as
$$
\sgn\of\mu\defeq\pas{-1}^{\numof{\text{crossings of }\mu}}.
$$
(See the right picture in \figref{fig:Kn-matchings}).

Now arrange the weights $a_{i,j}\defeq\weight\of{\setof{v_i,v_j}}$ % ($i<j$ by convention)
in an
upper triangular array $A=\pas{a_{i,j}}_{1\leq i < j \leq n}$:
%, $a_{i,j}\defeq\weight\of{\setof{i,j}}$.
The Pfaffian of $A$ is defined as
\begin{equation}
\label{eq:dfn-pfaffian}
\pfaffian\of{A}\defeq\sum_{\mu\in\matchings{\completeG V}}\sgn\of\mu\weight\of\mu,
\end{equation}
where the sum runs over all \pmatch s of $\completeG V$.

Since we \EM{always} view $\completeG V$ as weighted graph (with some weight function $\weight$),
we also
write $\pfaffian\of{\completeG V}$, or even simpler $\pfaffian\of{V}$, instead of
$\pfaffian\of{A}$. Moreover, since an upper triangular matrix $A$ determines uniquely
a skew symmetric matrix $A^\prime$ (by letting $A^\prime_{i,j}=A_{i,j}$ if $j>i$ and
$A^\prime_{i,j}=-A_{j,i}$ if $j<i$), we also write $\pfaffian\of{A^\prime}$ instead of
$\pfaffian\of{A}$. We set $\pfaffian\of{\emptyset}\defeq 1$ by definition.

%(Of course, this definition generalizes to the case of \EM{arbitrary} ordered vertex sets
%$\vertices\of{K_n}=\setof{v_1,\dots,v_{n}}$ in the obvious way.)
\end{dfn}

With regard to the identities for \pmatch s we are interested in, an edge not present
in graph $G$ may safely be added if it is given weight zero. Hence \EM{every} simple finite weighted graph $G$ may be viewed as a subgraph
(in general not an induced subgraph!) of $\completeG V$ with $\vertices\of{\completeG V}=\vertices\of G$, where
the weight of edge $e$ in $\completeG V$ is defined to be
\bit
\item $\weight\of e$ in $G$, if $e\in\edges\of G$,
\item zero, if $e\not\in\edges\of G$.
\eit
Keeping that in mind, we also write $\pfaffian\of{G}$ (or $\pfaffian\of{V}$, again)
instead of $\pfaffian\of{\completeG V}$.

\begin{figure}
\caption{Pfaffians according to Definition~\ref{dfn:Stembridge}: The left picture shows $\completeG 4$, drawn in the specific way described in
Definition~\ref{dfn:Stembridge}. The right picture shows the \pmatch\ 
$\mu=\setof{\setof{v_1,v_3},\setof{v_2,v_4}}$: Since there is precisely one crossing of
edges in the picture, $\sgn\of\mu = \pas{-1}^1 = -1$.}
\label{fig:Kn-matchings}
\input graphics/Kn-matchings
\end{figure}

The following simple observation is central for many of the following arguments.
\begin{obs}[contribution of a single edge to the sign of some \pmatch]
\label{obs:arcsnsigns}
Let $V=\pas{v_1,\dots,v_{2n}}$.
Removing an edge $e=\setof{v_i,v_j}$, $i<j$, together with the
vertices $v_i$ and $v_j$,
from some \pmatch\ $\mu$ of $\completeG{V}$ gives a \pmatch\ $\overline{\mu}$ of
the complete graph on the remaining vertices
$\pas{v_1,v_2,\dots,v_{2n}}\setminus\setof{v_i,v_j}$, and the change in sign
from $\mu$ to $\overline{\mu}$ is determined by the parity of
the number of vertices
\EM{lying between} $v_i$ and $v_j$ (see \figref{fig:arcsnsigns}).
%$\pas{-1}^{\numof{\text{vertices between $v_i$ and $v_j$}}}$.
By the ordering of the vertices, $\numof{\text{vertices between $v_i$ and $v_j$}}=j-i-1$, whence we have:
$$
\sgn\of{\mu} = \sgn\of{\overline{\mu}}\cdot\pas{-1}^{j-i-1}.
$$
\end{obs}

\begin{figure}
\caption{The contribution of edge $e=\setof{v_i,v_j}$ to the sign of the % perfect
matching $\pi$ amounts to $\pas{-1}^{\numof{\text{vertices between $v_i$ and $v_j$}}}$
(which is the same as $\pas{-1}^{i-j-1}$
if the vertices $\setof{v_1,v_2,\dots,v_{2n}}$ appear in ascending order).}
\label{fig:arcsnsigns}
\input graphics/observation
\end{figure}

\subsection{Cayley's Theorem and the long history of superposition of matchings}

The following Theorem of Cayley \cite{cayley:1849} is well--known. Stembridge
presented a beautiful proof (see \cite[Proposition 2.2]{stembridge:1990}) which was
based on superposition of matchings. Basically the same proof was already found
in the 19th century.
%
%The term \EM{superposition of
%matchings} for this combinatorial interpretation was coined by Knuth \cite{knuth:1996}.
%Apparently, it goes back to the 19th century.
We cite from \cite{knuth:1996}:
\begin{quote}
{\em 
An elegant graph-theoretic proof of Cayley's theorem%, somewhat
% analogous to the derivation of (1.0) above, 
 \dots was found by Veltmann in 1871
\cite{veltmann:1871}
and independently by Mertens in 1877
\cite{mertens:1877}.
Their proof anticipated 20th--century studies on the superposition of two
matchings, and the ideas have frequently been rediscovered.
}
\end{quote}

\begin{thm}[Cayley]
\label{thm:Cayley}
Given an upper triangular array $A=\pas{a_{i,j}}_{1\leq i < j \leq n}$, extend it
to a skew symmetric matrix $A^\prime=\pas{a^\prime_{i,j}}_{1\leq i,j \leq n}$ by setting
$$
a^\prime_{i,j} =
\begin{cases}
a_{i,j} & \text{ if } i < j, \\
-a_{i,j} & \text{ if } i > j, \\
0 & \text{ if } i = j.
\end{cases}
$$
Then we have:
\begin{equation}
\label{eq:pfaffiansquared-determinant}
\pas{\pfaffian\of{A}}^2 =
\det{\pas{A^\prime}}.
\end{equation}
\end{thm}

% Omit this in short version of the paper:
\iflongversion
\begin{proof}
By the definition of the determinant, we may view $\det{\pas{A^\prime}}$ as the generating
function of the symmetric group $\symm_n$
\begin{equation}
\label{eq:det-as-GF}
\det{\pas{A^\prime}} = \sum_{\pi\in\symm_n}\sgn\of\pi\weight\of\pi,
\end{equation}
where the weight $\weight\of\pi$ of a
permutation $\pi\in\symm_n$ is given as
$$
\weight\of\pi\defeq%\sgn\of{\pi}
\prod_{i=1}^n a^\prime_{i,\pi\of{i}}.
%\weight\of{v_i,v_{\pi\of i}}.
$$
The proof proceeds in two steps:

\step{1}
Denote by $\symmeven{n}$ the set of permutations $\pi\in\symm_n$ where the cycle decomposition
of $\pi$ \EM{does not} contain a cycle of odd length. Then we claim:
\begin{equation}
\label{eq:det-as-GF-even}
\det{\pas{A^\prime}} = \sum_{\pi\in\symmeven{n}}\sgn\of\pi\weight\of\pi.
\end{equation}
To prove this, we define a weight--preserving and sign--reversing involution
on the set of all permutations $\pi\in\pas{\symm_n\setminus\symmeven{n}}$ which
\EM{do} contain a cycle of odd length: % $\varphi$
%which preserves the absolute weight, but reverses the sign:
of all odd--length cycles in $\pi$, choose the one which contains the smallest
element $i$,
$$
\kappa_1=\pas{i,\pi\of{i},\pi^2\of{i},\dots,\pi^{2m}\of{i}},
$$
and replace it by its inverse
$$
\kappa_1^{-1}=\pas{\pi^{2m}\of{i},\pi^{2m-1}\of{i},\dots,\pi\of{i},i}.
$$
This operation obviously is an involution:
$$
\pi=\underbrace{\pas{\kappa_1,\kappa_2,\dots}}_{\text{cycle decomposition}}
\leftrightarrow\;
\pi^\prime=\underbrace{\pas{\kappa_1^{-1},\kappa_2,\dots}}_{\text{cycle decomposition}}.
$$
Since
$\weight\of{\pi^\prime}=-\weight\of\pi$ and $\sgn\of{\pi^\prime}=\sgn\of\pi$,
this involution is weight--preserving and sign--reversing.
So the terms corresponding to $\pas{\symm_n\setminus\symmeven{n}}$ in the \RHS\ of \eqref{eq:det-as-GF} sum up to zero, which proves \eqref{eq:det-as-GF-even}.

\step{2}
We shall construct a weight-- and sign--preserving bijection
between the terms 
\bit
\item $\sgn\of{\pi}\weight\of\pi$ corresponding to the determinant
as given in \eqref{eq:det-as-GF-even} (i.e., $\pi\in\symmeven{n}$)
\item
and %the terms
$\sgn\of\mu\weight\of\mu\sgn\of\nu\weight\of\nu$ corresponding to the square of the
Pfaffian in \eqref{eq:pfaffiansquared-determinant}.
\eit

To this end, consider the \EM{unique} cycle decomposition of $\pi$
\begin{equation}
\label{eq:pi-canonical}
\pi=
	\pas{i_1,\pi\of{i_1},\pi^2\of{i_1},\dots}%\concat
	\pas{i_2,\pi\of{i_2},\pi^2\of{i_2},\dots}%\concat
	\cdots
	\pas{i_m,\pi\of{i_m},\pi^2\of{i_m},\dots},
\end{equation}
where
\bit
\item $i_k$ is the smallest element in its cycle,
\item and
$i_1<i_2<\cdots<i_m$.
\eit
Visualize $\pi$ as directed graph as follows: represent
$$i_1,\pi\of{i_1},\dots,i_2,\pi\of{i_2},\dots$$ \EM{in this order}
(i.e., in the order in which the elements appear in \eqref{eq:pi-canonical}) by vertices
$$\pas{1,0},\pas{2,0},\dots,\pas{n,0}$$
in the plane. 
Call $\pas{1,0},\pas{3,0},\dots$ the
\EM{odd} vertices, and call $\pas{2,0},\pas{4,0},\dots$ the
\EM{even} vertices. Note that $\pi$ maps elements corresponding to
even vertices to elements corresponding to odd vertices, and
vice versa. If some element $i$ corresponds to an odd vertex $v$,
then draw a blue semicircle arc in the upper halfplane from $v$
to the even vertex $w$ which corresponds to $\pi\of i$.
If some element $j$ corresponds to an even vertex $s$,
then draw a red semicircle arc in the lower halfplane from $s$
to the odd vertex $t$ which corresponds to $\pi\of j$.
(See the left picture in \figref{fig:cayley} for an illustration.)

Note that if we forget the \EM{orientation} of the arcs, we simply
have a superposition $\pas{{\blue\mu},{\red\nu}}$ of a blue
and a red matching. Some of the arcs are co--oriented (i.e., they
point from left to right), and some are contra--oriented (i.e., they
point from right to left). Define the sign of any such \EM{oriented}
superposition of matchings by
\begin{equation}
\label{eq:sign-oriented}
\sgn\of{{\blue\mu},{\red\nu}}\defeq
\sgn\of{\blue\mu} \cdot
\sgn\of{\red\nu} \cdot
\pas{-1}^{
	%\numof{\text{crossings of }
	\numof{\text{contra--oriented arcs in }{\blue\mu}\cup{\red\nu}}
	}
\end{equation}
and observe that for the particular oriented superposition of matchings
obtained by visualizing permutation $\pi$ as above, this definition
gives precisely $\sgn\of\pi$.
(Again, see the left picture in \figref{fig:cayley} for an illustration.)

Furthermore, observe that
%since $\pi\in\symmeven{n}$, we may view the term $\sgn\of\pi\weight\of\pi$
%% corresponding to $\pi$ 
%as the signed weight of a
%superposition of matchings $\pas{\mu,\nu}$, where $\mu\cup\nu=\setof{\setof{1,\pi\of1},\dots,\setof{n,\pi\of n}}$. %(recall Observation~\ref{obs:interpretation}).
%More precisely,
with notation
$$\descof\pi\defeq\cardof{\setof{i:\;1\leq i\leq n,\,\pi\of{i}>i}}$$
we can rewrite the weight of $\pi$ as
\begin{equation}
\label{eq:weight-pi-rewritten}
\weight\of\pi=%\sgn\of\pi\cdot
\weight\of{\mu}\cdot\weight\of{\nu}
\pas{-1}^{\descof\pi}.
\end{equation}

However, the vertices do not appear in their original order. Clearly, we can obtain
the original ordering 
by interchanging neighbouring vertices $\pas{i,0}$ and $\pas{i+1,0}$ whose
corresponding elements appear in the wrong order, one after another, together with
the arcs being attached to them: see the right picture in
\figref{fig:cayley} and observe that this interchanging of vertices does not change
the sign as defined in \eqref{eq:sign-oriented}. Note that after finishing this ``sorting procedure'',
the number of contra--oriented arcs equals $\descof\pi$, so we have altogether
$$
\sgn\of\pi = \sgn\of{\blue\mu} \cdot
\sgn\of{\red\nu} \cdot
\pas{-1}^{\descof\pi}.
$$
Together with \eqref{eq:weight-pi-rewritten}, this amounts to
$$
\sgn\of\pi\cdot\weight\of\pi = 
	\pas{\sgn\of{\blue\mu}\cdot\weight\of{\blue\mu}} \cdot
	\pas{\sgn\of{\red\nu}\cdot\weight\of{\red\nu}},
$$
the \RHS\ of which obviously corresponds to a term in
$\pas{\pfaffian\of{A}}^2$.

On the other hand, every term in $\pas{\pfaffian\of{A}}^2$
corresponds to some superposition of matchings $S=\pas{{\blue\mu},{\red\nu}}$. For a
bicoloured cycle $C$ in $S$, identify the smallest vertex $v\in C$ and consider the
unique blue edge $\setof{v,w}$ in $C$. Orienting all bicoloured cycles $C$ such that this blue
edge points ``from $v$ to $w$'' gives an \EM{oriented} superposition of matchings, from
which we obtain a permutation without odd--length cycles and with the same weight 
and the same sign (by simply reversing the above ``sorting procedure'').
\end{proof}

\begin{figure}
\caption{Illustration for Cayley's Theorem. The left picture shows a cycle $c$ of length 8,
whose smallest element is $i$, i.e.,
$$c=\pas{i,\pi\of i, \pi^2\of{i},\dots,\pi^7\of i},$$
drawn as superposition of two directed matchings. Note that there is no crossing and
precisely one contra--oriented arc, whence, according to \eqref{eq:sign-oriented},
$c$ contributes $\pas{-1}$ to the sign of $\pi$,
as it should be for an even--length cycle. The right picture shows the effect of changing
the position of two neighbouring vertices $a$ and $b$. For \EM{both} matchings (red and blue;
blue arcs appear as dashed lines), we have:
%\begingroup
%\advance\leftskip by -3cm
%\bit
%\item
\newline$\bullet$
the number
of crossings changes by $\pm1$ if $a$ and $b$ belong to different arcs,
%\item
\newline$\bullet$
and if $a$ and $b$
belong to the same arc $e$, the orientation of $e$ is changed.
%\eit
%\endgroup
Since this
amounts to a change in sign for the red matching \EM{and} for the blue matching, the total
effect is that the sign \EM{does not change}.}
\label{fig:cayley}
\input graphics/cayley
\end{figure}

\fi

%It is common to define Pfaffians by \eqref{eq:dfn-pfaffian} for skew symmetric matrices
%(instead of upper triangular arrays). So by Cayley's Theorem, we have up to sign:
%$$
%\pfaffian\of{A} \defeq \sqrt{\det{\pas{A^\prime}}}.
%$$
% where, of course, the sign of the square root has to be determined properly.
\subsection{A corollary to Cayley's Theorem}
The following Corollary is an immediate consequence of Cayley's theorem. However,
we shall provide a direct ``graphical'' proof.
% (see \figref{fig:pfaffian-bipartite}).

Assume that the set of vertices $V$
is partitioned in two disjoint sets $A=\pas{a_1,\dots,a_m}$ and $B=\pas{b_1,\dots,b_n}$
such that the \EM{ordered} set $V$ appears as $\pas{a_1,\dots,a_m,b_1,\dots,b_n}$. Denote the complete
bipartite graph on $V$ (with set of edges $\setof{\setof{a_i,b_j}:1\leq i \leq m, 1\leq j \leq n}$) by $\completeG{A:B}$.
(For our purposes, we may view $\completeG{A:B}$ as the complete graph $\completeG{A\cup B}$,
where $\weight\of{\setof{a_i,a_j}}=\weight\of{\setof{b_k,b_l}}=0$ for all $1\leq i < j \leq m$
and for all $1\leq k < l \leq n$.) We introduce the notation
$$
\pfaffian\of{A,B}\defeq\pfaffian\of{\completeG{A:B}}.
$$

\begin{cor}
\label{cor:cayley}
Let $A=\pas{a_1,\dots,a_m}$ and $B=\pas{b_1,\dots,b_n}$ be two disjoint ordered sets.
% be a \EM{bipartite} graph:  $\vertices\of G=\setof{a_1,\dots,a_m,b_1,\dots,b_n}$
%and each $e\in\edges\of G$ is of the form $\setof{a_i,b_j}$.
Then we have
%$\vertices\of G=A\cup B$, $A\cap B=\emptyset$ and
%$e\in\edges\of G\implies\pas{e\cap A\neq\emptyset \text{ and } e\cap B\neq\emptyset}$.
$$
\pfaffian\of{A,B}=%\pfaffian\of{\completeG{A:B}}=
\begin{cases}
\pas{-1}^{\binom{n}{2}}\det\of{\weight\of{a_i,b_j}}_{i,j=1}^{n} &\text{ if } m=n,\\
0 & \text{ otherwise.}
\end{cases}
$$
\end{cor}
\begin{proof}
If $m\neq n$, then there is no \pmatch\ in $\completeG{A:B}$, and thus the Pfaffian clearly is $0$.

If $m=n$,
consider the $n\times n$--matrix $M\defeq\pas{\weight\of{\setof{a_i,b_{n-j+1}}}}_{i,j=1}^{n}$.
Note that for every permutation $\pi\in\symm_n$, the corresponding term in the expansion of
$\det\of M$ may be viewed as the signed weight of a certain \pmatch\ $\mu$ of $\completeG{A:B}$
(see
\figref{fig:pfaffian-bipartite}). Recall that $\sgn\of\pi=\pas{-1}^{\inv\of\pi}$, where
$\inv\of\pi$ denotes the number of inversions of $\pi$, and observe that \EM{inversions}
of $\pi$ are in one--to--one--correspondence with \EM{crossings} of $\mu$. Thus
$$
\pfaffian\of{\completeG{A:B}}=\det\of{M},
$$
and the assertion follows by reversing the order of the columns of $M$.
\end{proof}

\begin{figure}
\caption{Illustration for Corollary~\ref{cor:cayley}: Consider the $4\times 4$--matrix
$M\defeq\pas{\weight\of{\setof{a_i,b_{5-j}}}}_{i,j=1}^{4}$ and the permutation
$\pi=\pas{2,3,4,1}$ in $\symm_4$. The left pictures shows $\pi$ as (bijective) function
mapping the set $\setof{1,2,3,4}$ onto itself: It is obvious that the ``arrows''
indicating the function constitute a \pmatch\ $\mu$. The right picture shows
the same matching $\mu$ drawn in the specific way of Definition~\ref{dfn:Stembridge}.
Inversions of $\pi$ are in one--to--one correspondence with crossings of $\mu$, whence we see:
$$\sgn\of\pi=\sgn\of\mu.$$}
\label{fig:pfaffian-bipartite}
\input graphics/pfaffian-bipartite
\end{figure}

% Omit this in short version of the paper:

\subsection{A generalization of Corollary~\ref{cor:cayley}}
We may use Observation~\ref{obs:arcsnsigns} to prove another identity involving
Pfaffians. To state it conveniently, we need some further notation.

Assume that the ordered set of vertices $V$
appears as $\pas{a_1,\dots,a_m,b_1,\dots,b_n}$ for disjoint sets $A=\pas{a_1,\dots,a_m}$
and $B=\pas{b_1,\dots, b_n}$. Consider the complete graph $\completeG V$ and delete (or assign weight zero to) all edges
joining two vertices from $A$. Call the resulting graph the \EM{semi--bipartite}
graph $\semibip{A}{B}$. (Note that every \pmatch\ $\mu$ in $\semibip{A}{B}$
constitutes an \EM{injective} mapping $A\to B$.)

Let $M=\pas{m_1,m_2,\dots,m_n}$ be some ordered set. For some arbitrary subset 
$X=\pas{m_{i_1},\dots,m_{i_k}}\subseteq M$, denote by $\SETSUM{X}{M}$ the sum of the
\EM{indices} $i_j$ of $X$:
%of indices $1\leq{i_1}<{i_2}<\dots<{i_k}\leq n$, denote
$$\SETSUM{X}{M}\defeq i_1 + i_2 + \cdots + i_k.$$
(Recall that
subsets ``inherit'' the ordering in our presentation, i.e., $i_1<i_2<\dots<i_k$:
we might also call $X$ a \EM{subword} of $M$.)

% Omit this in short version of the paper:
\iflongversion
\begin{rem}
Note that this concept is related to % sign of
the permutation $\pi\in\symm_M$ which moves $m_{i_j}$ to
position $j$,
$j=1\dots,k$, while leaving the order of all remaining elements
in $M\setminus X$ unchanged:
$$
\sgn\of\pi =\pas{-1}^{{i_1}-1+{i_2}-2+\dots+{i_k}-k}=\pas{-1}^{\SETSUM{X}{M}-\binom{k+1}{2}}.
$$
(This observation provides the translation to the presentation of Pfaffians
given in \cite{knuth:1996} and \cite{hamel:2001}, it is not needed for our presentation.)
\end{rem}
\fi

\begin{cor}
\label{cor:semi-bipartite}
Let $V=\pas{a_1,\dots,a_m,b_1,\dots,b_n}$, $A=\pas{a_1,\dots,a_m}$ and $B=\pas{b_1,\dots, b_n}$.
For every subset $Y=\setof{b_{k_1},b_{k_2},\dots,b_{k_m}} \subseteq B$,
denote by $M_Y$ the $m\times m$--matrix
$$M_Y\defeq\pas{\weight\of{\setof{a_i,b_{k_j}}}}_{i,j=1}^m.$$
Then we have 
\begin{equation}
\label{eq:semi-bipartite}
\pfaffian\of{\semibip{A}{B}}=
% \pas{-1}^{\binom{m}{2}}
\pas{-1}^{m}
\sum_{\substack{Y\subseteq B,\\ \cardof{Y}=m}}
\pas{-1}^{\SETSUM{Y}{B}}\cdot\pfaffian\of{B\setminus Y}\cdot\det\of{M_Y}.
\end{equation} 
\end{cor}
\begin{proof}
For every \pmatch\ $\rho$ in $\semibip{A}{B}$, let $Y\subseteq B$ be the set
of vertices which are joined with a vertex in $A$ by some edge in $\rho$. Note
that $\cardof{Y}=\cardof{A}=m$ (if such \pmatch\ exists), and observe that
$\rho$ may be viewed as a superposition of \pmatch s, namely
\bit
\item a red matching ${\red\mu}$ in the complete bipartite graph $\completeG{A,Y}$
\item and a blue matching ${\blue\nu}$ in the complete graph $\completeG{B\setminus Y}$,
\eit
where
$$
\weight\of\rho = \weight\of{\red\mu}\cdot\weight\of{\blue\nu}.
$$
For an illustration, see \figref{fig:cor-cayley}. Note that the crossings of $\rho$
are partitioned in
\bit
\item crossings of two edges from $\red\mu$,
\item crossings of two edges from $\blue\nu$
\item and crossings of an edge from $\red\mu$ with an edge from $\blue\nu$,
\eit
whence we have
$$\sgn\of\rho =
\pas{-1}^{\numof{\text{crossings of an edge from $\red\mu$ \EM{and} an edge from $\blue\nu$}}}
\cdot\sgn\of{{\red\mu}}\cdot\sgn\of{{\blue\nu}}.$$

\goodbreak
Assume that $Y=\pas{b_{k_1},b_{k_2},\dots,b_{k_m}}$ and observe that \EM{modulo 2}
the number of crossings
\bit
\item of the edge from $\red\mu$ which ends in $b_{k_j}$
\item with edges from $\blue\nu$
\eit equals the number of vertices of $B\setminus Y$ which lie to the \EM{left} of $b_{k_j}$,
which is $k_j - j$. Hence we have
$$\sgn\of\rho\cdot\sgn\of{{\red\mu}}\cdot\sgn\of{{\blue\nu}} =
\pas{-1}^{(k_1-1)+(k_2-2)+\dots+(k_m-m)} =
\pas{-1}^{\SETSUM{Y}{B} -\binom{m+1}{2}}.$$
From this we obtain
\begin{equation}
\label{eq:proof-semibip}
\pfaffian\of{\semibip{A}{B}}=
\pas{-1}^{m}
\sum_{\substack{Y\subset B,\\ \cardof{Y}=m}}
\pas{-1}^{\SETSUM{Y}{B}}\cdot
\pfaffian\of{B\setminus Y}\cdot
\pas{-1}^{\binom{m}{2}}\cdot
\pfaffian\of{A,Y},
\end{equation}
which by Corollary~\ref{cor:cayley} equals \eqref{eq:semi-bipartite}.
\end{proof}

\begin{figure}
\caption{Illustration for Corollary~\ref{cor:semi-bipartite}. Consider the ordered set
of vertices 
$$V=\setof{a_1,a_2,a_3;b_1,b_2,\dots,b_7}.$$
The picture shows the matching
$$
\rho=
\setof{\setof{a_1,b_5},\setof{a_2,b_2},\setof{a_3,b_7},\setof{b_1,b_4},\setof{b_3,b_6}}
$$
in $\semibip{A}{B}$, where $A=\setof{a_1,a_2,a_3}$ and $B=\setof{b_1,b_2,\dots,b_7}$.
Let $Y=\setof{b_2,b_5,b_7}$ and observe that $\rho$ may be viewed as superposition of
% \bit
% \item
the
red matching ${\red\mu}=\setof{\setof{a_1,b_5},\setof{a_2,b_2},\setof{a_3,b_7}}$ in
$\completeG{A,Y}$
% \item
and the blue matching ${\blue\nu}=\setof{\setof{b_1,b_4},\setof{b_3,b_6}}$ in
$\completeG{B\setminus Y}$ (blue edges are drawn as dashed lines).
% \eit
All crossings in $\rho$
are indicated by circles; the crossings which are \EM{not} present in $\red\mu$ or in $\blue\nu$
are indicated by two concentric circles.
}
\label{fig:cor-cayley}
\input graphics/cor-cayley
\end{figure}

\section{\Pmatch s and Pfaffians: The Kasteleyn--Percus method}
\label{sec:kasteleyn-percus}

We now present the Kasteleyn--Percus method for
the enumeration of matchings in plane graphs, following closely (and
refining slightly) the
exposition in \cite{kasteleyn:1967}.

The main idea is simple: If we \EM{disregard the signs} of the terms,
the Pfaffian $\pfaffian\of{D\of{G}}$, by definition, encompasses
the same terms as the generating function $\nofmatchings{G}$ of matchings in
$G$. If it is possible to modify the weight function $\weight$ by
introducing signs such that \EM{for all} matchings $\mu$ in $G$ the modified weight
function $\weight^\prime$ is
$$
\weight^\prime\of\mu = \sgn\of\mu\weight\of\mu,
$$
then the Pfaffian (\EM{for the modified weight--function} $\weight^\prime$) would be equal to $\nofmatchings{G}$
(\EM{for the original weight--function} $\weight$).

So let $G$ be some \loopless\ graph with weight function $\weight$ and assume
some \EM{orientation} $\orientation$ on the pairs
of vertices of $G$:% in the following sense:
\begin{equation}
\orientation: \vertices\of{G}\times\vertices\of{G}\to \setof{1,-1}
\text{ such that } \orientation\of{v,u} = - \orientation\of{u,v}.
\end{equation}  
Consider the skew--symmetric square matrix $D\of{G,\orientation}$
with row and column
indices corresponding to the vertices $\setof{v_1,\dots,v_n}$ of $G$
(in some arbitrary
order) and entries
\begin{align*}
d_{i,i} &=0,\\
d_{i,j} & =\orientation\of{v_i,v_j}\times
			\sum_{\substack{e\in\edges\of{G},\\e=\setof{v_i,v_j}}}\weight\of e
			\text{ for }i\neq j.
\end{align*}

%If we \EM{disregard the signs} of the terms, the Pfaffian
%$\pfaffian\of{D\of{G,\orientation}}$, by definition, encompasses
%the same terms as the generating function $\nofmatchings{G}$ of matchings in
%$G$.
Clearly, the weights $\weight^\prime\of\mu$ of the terms in the Pfaffian $\pfaffian\of{D\of{G,\orientation}}$ differ from the weights $\weight\of\mu$
of the terms in the Pfaffian $\pfaffian\of{D\of{G}}$ (i.e., $G$ without orientation)
only by a \EM{sign} which depends on the orientation
$\orientation$. So if we find an %Kasteleyn's idea was to look for an
orientation $\orientation$
of $G$ under which \EM{all} the terms in the Pfaffian have the \EM{same} sign
$\pas{-1}^m$, i.e., for all matchings $\mu$ of $G$ we have
$$
\sgn\of\mu\weight^\prime\of\mu = \pas{-1}^m\weight\of\mu,
$$
% If such an orientation exists,
then the generating function $\nofmatchings{G}$
% of the matchings in graph $G$
is equal to %the Pfaffian
$\pas{-1}^m\pfaffian\of{D\of{G,\orientation}}$.

All terms in the Pfaffian $\pfaffian\of{D\of{G,\orientation}}$
have the \EM{same} sign if and only if all the terms in the squared Pfaffian
$\pfaffian\of{D\of{G,\orientation}}^2$ have the \EM{positive} sign, which
by Cayley's theorem (stated previously as Theorem~\ref{thm:Cayley})
is equivalent to all the terms in the determinant $\det\of{D\of{G,\orientation}}$ being positive. According to the proof of Cayley's theorem, the terms in this determinant
correspond to permutations $\pi$ with cycle decompositions where every cycle
has \EM{even} length. Since an even--length cycle contributes the
factor $(-1)$ to the sign of $\pi$, i.e.,
$$
\sgn\of\pi=\pas{-1}^{\text{number of even--length cycles in $\pi$}},
$$
the overall sign of the term in the determinant certainly will be positive
if \EM{each} of the even--length cycles contributes an offsetting factor
$(-1)$, i.e., if every even--length cycle in $\pi$ contains an \EM{odd} number
of elements $d_{i,\pi\of i}$ with negative sign. Note that this condition
is always fulfilled for cycles of length 2: Exactly one of the elements
in $\pas{d_{i,j} d_{j,i}}$ has the negative sign.

These considerations can be restated in terms of the graph $G$.
\begin{dfn}
Let $G$ be some \loopless\ graph with weight function $\weight$ and \EM{orientation} $\orientation$.
% Recall that %, according to Observation~\ref{obs:interpretation},
The superposition of two arbitrary \pmatch s of $G$
yields a covering of the bicoloured graph
$B=\bicoloured{G}{\thereds}{\theblues}$, $\thereds=\theblues=\emptyset$
(i.e., there are no coloured vertices), with % double edges and
even--length cycles. (Recall that a superposition of matchings in $G$
corresponds to a term in the squared Pfaffian
$\pfaffian\of{D\of{G,\orientation}}^2$, and the corresponding covering with
even--length cycles corresponds to a term in the determinant 
$\det\of{D\of{G,\orientation}}$; see the proof of Cayley's theorem.)
%Consider the ``inherited'' orientation $\orientation$ for $B$ (i.e.,
%the same orientation as in $G$.)

A cycle in $B$ of length $>2$, which arises from the
superposition of two \pmatch s of $G$, corresponds to a ``normal''
even--length cycle $C$ in $G$: We call such cycle $C$ a
\EM{superposition cycle}.
%Note that a superposition cycle of length $>2$ , while a superposition cycle of length $2$
%corresponds to a ``doubled edge'' in $G$.

The orientation $\orientation$ is called
\EM{admissible} if in \EM{every} superposition cycle $C$ there is an \EM{odd}
number of edges \EM{co--oriented} with $\orientation$ (and an \EM{odd}
number of edges \EM{contra--oriented} with $\orientation$; here, co-- and
contra--orientation refer to some
arbitrary but fixed orientation of the even--length cycle $C$).
\end{dfn}

These considerations can be summarized as follows \cite[Theorem~$\brk{1}$ on page 92]{kasteleyn:1967}:
\begin{thm}[Kasteleyn]
\label{thm:kasteleyn1}
Let $G$ be a \loopless\ graph with weight function $\weight$. If $G$
has an \EM{admissible} orientation
$\orientation$, then the total weight of all \pmatch s
of $G$ equals the Pfaffian of $D\of{G,\orientation}$ up to sign:
\begin{equation}
\nofmatchings{G} = \pm\pfaffian\of{D\of{G,\orientation}}.
\end{equation}
\end{thm}

\subsection{Admissible orientations for planar graphs}
While the existence of an admissible orientation is not guaranteed
in general, for a \EM{planar graphs} $G$ such orientation can be constructed
\cite{kasteleyn:1967}.

\begin{figure}
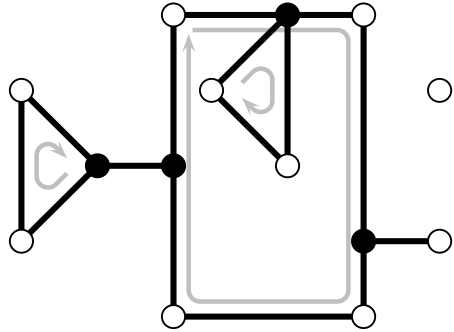

\caption{Decomposition of a graph in $2$--connected blocks, bridges and
isolated vertices: The cut--vertices
are drawn as black circles. The graph shown here is planar, each of its
three $2$--connected blocks consists of a single cycle. The \EM{clockwise} orientation
of these cycles (in the given embedding) is indicated by grey arrows.
}
\label{fig:block-decomposition}
\input graphics/block-decomposition
\end{figure}

For this construction, we need some facts from graph theory.
Let $G$ be a graph. If there are two
different vertices $p\neq q\in\vertices\of G$ belonging to the same connected
component $G^\star$ of $G$, such that there is \EM{no} cycle in $G^\star$ that
contains \EM{both} vertices $p$ and $q$ (i.e., $G^\star$ is \EM{not}
$2$--connected), then by Menger's Theorem (see, e.g.,
\cite[Theorem 3.3.1]{diestel:2000}) there exists a vertex $v$ in $G^\star$
such that $\subgraph{G^\star}{\setof{v}}$ is disconnected: Such vertex
$v$ is called an \EM{articulation vertex} or \EM{cutvertex}. The whole
graph $G$ is subdivided by its cutvertices, in the following
sense: Each cutvertex connects two or more \EM{blocks}, i.e.,
maximal connected subgraphs that \EM{do not contain} a cutvertex.
Such blocks are
\bit
\item either maximal $2$--connected subgraphs $H$ of $G$ (i.e., for every
	pair of different
	vertices $p,q\in\vertices\of{H}$ there \EM{exists} a cycle in $H$ that
	contains both $p$ and $q$),
\item or single edges (called \EM{bridges}),
\item or isolated vertices.
\eit
See \figref{fig:block-decomposition} for an illustration.

Since we deal with \EM{cycles} here, we are mainly interested in the
the $2$--connected blocks which are not isolated
vertices. (Clearly, a graph with an isolated vertex has no matching.)

In the following, assume that $G$ is planar and consider an arbitrary but
\EM{fixed} embedding of $G$
in the plane: So when we speak of $G$, we always mean ``$G$ in its
fixed planar embedding''.

For every $2$--connected block $H$ of $G$, consider the embedding ``inherited''
from $G$. Note that the boundary of a face of a $2$--connected planar graph
always is a \EM{cycle}. 

% We shall never consider the \EM{unbounded} face of the graph embedding here.
A cycle which is the boundary of a \EM{bounded} face (we shall never consider
the \EM{unbounded} face of the graph embedding here) of a
$2$--connected block $H$ of $G$ (looking at $H$ \EM{alone}, forgetting
the rest of $G$) is called a \EM{contour cycle} in $G$. The vertices of $G$
lying in the \EM{interior} of some % contour
cycle $C$  (i.e., lying
in the bounded region confined by $C$ in the fixed embedding) are called the
\EM{interior vertices} of $C$. For the number of
interior vertices of $C$ we introduce the notation $\nofinterior{C}$.

In the plane, consider the \EM{clockwise orientation}: This determines
a unique orientation for every % contour
cycle of $G$ (by choosing a ``center''
in the bounded region confined by $C$ in the fixed embedding of $G$
and traversing the edges of $C$ in the clockwise
orientation around this center). Given some arbitrary orientation $\orientation$ of the edges
of $G$, we call the edges of $C$ for which $\orientation$ coincides with the
clockwise orientation the \EM{co--oriented edges} of $C$ (the other edges
are called the \EM{contra--oriented edges} of $C$). For the number of
co--oriented edges of $C$ we introduce the notation $\nofcooriented{C}$.

Then we have \cite[Lemma~$\brk{2a}$ on page 93]{kasteleyn:1967}:
\begin{lem}
\label{lem:2a}
For each finite planar graph $G$ there is an orientation $\orientation$ of the
edges of $G$
such that for every \EM{contour} cycle $C$ of $G$
$$
\nofcooriented{C} \not\equiv \nofinterior{C} \pmod 2.
$$
%\bit
%\item the number of co--oriented edges of $C$
%\item and the number of interior vertices of $C$
%\eit
%are of opposite parity.
\end{lem}

\begin{proof}
For bridges (edges not belonging to a $2$--connected block) in $G$,
we may choose an arbitrary orientation.

For the remaining edges, we may construct the orientation by
considering \EM{independently} the $2$--connected blocks $H$ of $G$, one after
another. 

So let $H$ be a $2$--connected block with the embedding inherited from $G$.
For the moment, look at $H$ \EM{alone}, i.e., forget the rest of $G$. The algorithm is as follows:

Start by choosing an arbitrary contour
cycle $C_1$ in $H$ and choose an orientation for its edges such that the
number of co--oriented edges of $C_1$ and the number of interior vertices
of $C_1$ are of opposite parity (clearly, this is possible). Note
that in the inherited embedding of $H$, the union of $C_1$ and the face bounded by $C_1$ is a
\EM{simply connected region} in the plane (homeomorphic to the closed
disk).

Now repeat the following step until all edges of $H$ are oriented:
Assume that the edges belonging to contour cycles $C_1,C_2,\dots,C_n$
have already been oriented by our algorithm \EM{and} that the union
of all these cycles and corresponding faces is a
\EM{simply connected region} in the plane. Now choose a contour cycle
$C_{n+1}$ such that the union of the cycles $C_1,C_2,\dots,C_{n+1}$
together with their corresponding faces is again a 
\EM{simply connected region} in the plane. A moment's thought shows that
this is possible, and that the set of edges of $C_{n+1}$ which
were not yet oriented is non--empty: Clearly, for these remaining
edges we can choose an orientation such that the
number of co--oriented edges of $C_{n+1}$ and the number of interior vertices
of $C_{n+1}$ are of opposite parity.
\end{proof}

It is clear that \EM{every} cycle in the planar graph $G$ encircles one
or more faces. We use this simple fact to give the following generalization
of Lemma~\ref{lem:2a} \cite[Lemma~$\brk{2b}$ on page 93]{kasteleyn:1967}:
\begin{pro}
\label{prop:2b}
Let $G$ be a finite planar graph.
For the orientation constructed in Lemma~\ref{lem:2a} there holds:
For \EM{every} cycle $C$ of $G$ (not only contour cycles!) we have
$$
\nofcooriented{C} \not\equiv \nofinterior{C} \pmod 2.
$$
\end{pro}
\begin{proof}
We prove this by induction on the number $n$ of faces encircled by
the cycle $C$: For $n=1$, the assertion is true by Lemma~\ref{lem:2a}.

So assume the assertion to be true for all cycles encircling $n$ faces,
and consider some cycle $C$ encircling $n+1$ faces. Select one of these
faces $C_i$ such that the union of all other faces $\bigcup_{j=1,j\neq i}^{n+1} C_j$
(together with their corresponding contour cycles)
is also a simply connected region: A moment's thought shows that this is
always possible.

For the cycle $C^\prime$ encircling $\bigcup_{j=1,j\neq i}^{n+1} C_j$
and for the cycle $C_i$, the assertion is true by induction.
By construction, the edges belonging to \EM{both} $C^\prime$ \EM{and}
$C_i$ form a path %$p$
of length $k>0$,
$$
\pas{v_0,v_1,\dots,v_k},
$$
and the \EM{new} interior vertices of $C$ (which are not also interior
vertices of $C_i$ or $C^\prime$) are precisely $v_1,\dots v_{k-1}$.
Now observe that for every edge $e$ of $p$ we have: If $e$ is co--oriented in
$C_i$, then it is contra--oriented in $C^\prime$, and vice versa. Hence
we have
\begin{align*}
\nofcooriented{C} &=
\nofcooriented{C^\prime} +
\nofcooriented{C_i} -
k, \\
\nofinterior{C} &=
\nofinterior{C^\prime} +
\nofinterior{C_i} +
k-1.
\end{align*}
This proves the assertion.
\end{proof}

Now we obtain immediately the following result \cite[Theorem~$\brk{2}$ on page 94]{kasteleyn:1967}:
\begin{thm}[Kasteleyn]
\label{thm:kasteleyn2}
The orientation constructed in Lemma~\ref{lem:2a} is admissible: I.e.,
for every finite \EM{planar} graph $G$ there exists an admissible orientation.
%, such that for
%\EM{every} even--length cycle in $G$ the number of co--oriented edges
%is odd. (In particular, the orientation is \EM{admissible}).
\end{thm}
\begin{proof}
Simply observe that for a \EM{superposition cycle} $C$ of a \EM{planar} graph
the number of interior points is necessarily \EM{even}. (Recall that different
superposition cycles cannot have a vertex in common.) So the assertion follows
from Proposition~\ref{prop:2b}.
\end{proof}

Observe that the admissibility of the orientation is ``inherited'' by
certain induced subgraphs:
\begin{cor}
\label{cor:special-subgraphs}
Let $G$ be a finite planar graph with the admissible orientation $\orientation$
constructed in Lemma~\ref{lem:2a}. Let $C$ be some contour cycle of
$G$, and let $S=\setof{v_1,\dots,v_{2k}}$ be some set of $2k$ vertices of
$C$. Consider $H=\subgraph{G}{S}$ with the orientation $\orientation$
inherited from $G$: Then $\orientation$ is admissible for $H$.
\end{cor}
\begin{proof}
Simply note that the face bounded by $C$ in $G$ belongs to a bigger
face in $H$: The contour cycle for this face in $H$ is a cycle $C^\prime$ 
in $G$, for which $\nofcooriented{C^\prime}\not\equiv\nofinterior{C^\prime}\pmod 2$
holds in $G$ \EM{and} in $H$ (since the number of interior points of $C^\prime$
is decreased by $2k$ in $H$).
\end{proof}
So it seems that Proposition~\ref{prop:kuo} gives an identity for \EM{special}
Pfaffians which correspond to planar graphs $G$: Let $\orientation$ be
an admissible orientation for $G$, and assume the
same (inherited) orientation for induced subgraphs of $G$, then
\eqref{eq:kuos-proposition} translates to
\begin{multline}
\label{eq:kuos-pfaffian}
\pm\pfaffian\of{D\of{G,\orientation}}\cdot
\pfaffian\of{D\of{\subgraph{G}{\setof{a,b,c,d}},\orientation}}\\\pm
\pfaffian\of{D\of{\subgraph{G}{\setof{a,c}},\orientation}}\cdot
\pfaffian\of{D\of{\subgraph{G}{\setof{b,d}},\orientation}} =
\\
\pm\pfaffian\of{D\of{\subgraph{G}{\setof{a,b}},\orientation}}\cdot
\pfaffian\of{D\of{\subgraph{G}{\setof{c,d}},\orientation}}\\
\pm 
\pfaffian\of{D\of{\subgraph{G}{\setof{a,d}},\orientation}}\cdot
\pfaffian\of{D\of{\subgraph{G}{\setof{b,c}},\orientation}}
\end{multline}
for the ``proper'' choice of signs, by Theorem~\ref{thm:kasteleyn1}.

But it turns out that the ``proper'' choice of signs is ``always $+$''
or ``always $-$'', and that \eqref{eq:kuos-pfaffian} with this
choice of signs is in fact
an identity for Pfaffians \EM{in  general}, namely the special case
\cite[Equation (1.1)]{knuth:1996} of an identity \cite[Equation (1.0)]{knuth:1996} due to Tanner \cite{tanner:1878}. This can be made precise as follows:
\begin{dfn}
\label{dfn:wpit}
Assume an identity for Pfaffians of the form
$$
\sum_{i=1}^m\pfaffian\of{G_i}\cdot\pfaffian\of{H_i} =
\sum_{j=1}^n\pfaffian\of{G^\prime_j}\cdot\pfaffian\of{H^\prime_j},
$$
where all the graphs $G_i$, $H_i$, $G^\prime_j$ and $H^\prime_j$
involved are induced subgraphs of some ``supergraph'' $G$, with
$$
\vertices\of{G_i}\cup\vertices\of{H_i} =
\vertices\of{G^\prime_j}\cup\vertices\of{H^\prime_j} =
\vertices\of{G}
$$
for $i=1,\dots,m$ and $j=1,\dots, n$.
If this identity comes, in fact, from a \EM{sign-- and weight--preserving involution} which maps
% mapping
\bit
\item the family of superpositions of matchings corresponding to the
\EM{left--hand side}
\item to the family of superpositions of matchings corresponding to the
	\EM{right--hand side},
\eit then we say that the
identity is of the \EM{\wpit}.
\end{dfn}

\begin{rem}
An identity of the \wpit\ could, of course, be written in the form
$$
\sum_{i=1}^m\pfaffian\of{G_i}\cdot\pfaffian\of{H_i} -
\sum_{j=1}^n\pfaffian\of{G^\prime_j}\cdot\pfaffian\of{H^\prime_j} = 0,
$$
and if we view the left--hand side as the generating function of signed
weights of superpositions of matchings
$$\pas{-1}^d\cdot\sgn\of\mu\cdot\sgn\of\nu\cdot\weight\of\mu\cdot\weight\of\nu,$$
where
\bit
\item $d=0$ for the objects corresponding
to the \EM{unprimed} pairs $\pas{G_i,H_i}$
\item and
$d=1$ for the objects corresponding
to the \EM{primed} pairs $\pas{G^\prime_j,H^\prime_j}$,
\eit
then the involution according to Definition~\ref{dfn:wpit} would appear as
\EM{sign--reversing}.
\end{rem}

\begin{lem}% [Kasteleyn--Percus]
\label{lem:kasteleyn-percus}
If an identity for Pfaffians
%$$
%\sum_{i=1}^m\pfaffian\of{G_i}\cdot\pfaffian\of{H_i} =
%\sum_{j=1}^n\pfaffian\of{G^\prime_j}\cdot\pfaffian\of{H^\prime_j},
%$$
is of the \wpit, and if this identity
can be specialized in a way such that
\bit
\item the ``supergraph'' $G$ is planar with admissible orientation $\orientation$,
\item and the  inherited orientation $\orientation$ is also admissible for \EM{all}
	the induced subgraphs $G_i$, $H_i$ ($i=1,\dots, m$) and $G^\prime_j$,
	$H^\prime_j$ ($j=1,\dots, n$) of $G$,
\eit then this identity ``translates'' immediately
to the corresponding identity for matchings, i.e., to
$$
\sum_{i=1}^m\nofmatchings{G_i}\cdot\nofmatchings{H_i} =
\sum_{j=1}^n\nofmatchings{G^\prime_j}\cdot\nofmatchings{H^\prime_j}.
$$
\end{lem}
\begin{proof}
Construct the graph $J$ with vertex set
$$
\underbrace{\setof{\pas{G_1,H_1},\pas{G_2,H_2},\dots,\pas{G_m,H_m}}}_{\text{unprimed vertices}}
\disjuni
\underbrace{\setof{\pas{G^\prime_1,H^\prime_1},\pas{G^\prime_2,H^\prime_2},\dots,\pas{G^\prime_n,H^\prime_n}}}_{\text{primed vertices}}.
$$
Connect two vertices $\pas{G_i,H_i}$ and $\pas{G^\prime_j,H^\prime_j}$ of $J$
by an edge if and only if the involution maps some superposition of matchings
from $\pas{G_i,H_i}$ to some superposition of matchings from $\pas{G^\prime_j,H^\prime_j}$. Clearly, $J$ is bipartite; the bipartition is given by the
sets of primed and unprimed vertices.

Recall that all the superpositions of matchings
in the identity which correspond to some \EM{fixed} vertex of $J$ have the same sign, since
the inherited orientations are admissible.
Let $C$ be an arbitrary connected component of $J$:
It is easy to see that \EM{all} the terms in the identity corresponding to vertices
of $C$ have the same sign. So cancelling the sign, we see that the total weight
(without sign!)
of the unprimed vertices of $C$ equals the total weight (without sign!)
of the primed vertices
of $C$, i.e., $C$ corresponds to an identity for matchings.

This consideration can be applied to \EM{all} connected components of $J$, and the
sum of the corresponding identities for matchings gives the desired translation
of the Pfaffian identity of the \wpit.
% proves the assertion.
\end{proof}

In the following, we shall say that an identity for matchings 
of planar graphs \EM{follows}
from an identity for Pfaffians \EM{by the Kasteleyn--Percus method}, if it
can be obtained by the ``translation'' in the sense of Lemma~\ref{lem:kasteleyn-percus}.

\section{Overlapping Pfaffian identities}% due to Tanner, Ohta and Krattenthaler}
\label{sec:tanner}
We consider Pfaffian identities here which are all of the \wpit: This fact
will be obvious immediately from the proofs we give, since they are all based on
weight--preserving involutions.

The following theorem is due to H.W.L.\ Tanner
\cite{tanner:1878} (see \cite[Equation (1.0)]{knuth:1996}.
%We present it here in
%our ``language'':
\begin{thm}
\label{thm:tanner}
Let $\word\alpha=\pas{\alpha_1,\dots,\alpha_n}$ and $\word\beta=\pas{\beta_1,\dots,\beta_m}$ be %words over some alphabet $X$.
subsets of the ordered set $\gamma=\pas{\alpha_1,\dots,\alpha_n,\,\beta_1,\dots,\beta_m}$.
Let $k$, $1\leq k\leq m$, be arbitrary but fixed. Then there holds
\begin{equation}
\label{eq:tanner}
%\symbolf\brk{\alpha}\symbolf\brk{\alpha\beta} =
%	\sum_{y}s\of{\beta; x  y}
%		\symbolf\brk{\alpha x y}
%		\symbolf\brk{\alpha\beta\setminus x y}.
\pfaffian\of{\alpha}\pfaffian\of{\alpha\cup\beta} = \pas{-1}^k\cdot
	\sum_{\substack{j=1,\\j\neq k}}^m%s\of{\beta; x  y}
		%\pas{-1}^{\SETSUM{\setof{x,y}}{\beta}}
		\pas{-1}^{j-1}
		\pfaffian\of{\alpha\cup\setof{\beta_k, \beta_j}}
		\pfaffian\of{\pas{\alpha\cup\beta}\setminus\setof{\beta_k, \beta_j}}.
%	\text{ for all } x\in\beta.
\end{equation}
(Recall that all subsets \EM{inherit} the order of $\word\gamma$.)

This identity is of the \wpit.
\end{thm}

Observe that for the special case $\beta = \pas{a, b, c, d}$, \eqref{eq:tanner}
reads
\begin{multline*}
\pfaffian\of{\alpha}\cdot\pfaffian\of{\alpha\cup\pas{a, b, c, d}}
	   + \pfaffian\of{\alpha\cup\pas{a, c}}\cdot
		\pfaffian\of{\alpha\cup\pas{b, d}}
		 =\\
		\pfaffian\of{\alpha\cup\pas{a, b}}\cdot
		\pfaffian\of{\alpha\cup\pas{c, d}} +
		\pfaffian\of{\alpha\cup\pas{a, d}}\cdot
		\pfaffian\of{\alpha\cup\pas{b, c}},
\end{multline*}
which implies \eqref{eq:kuos-pfaffian} by the
Kasteleyn--Percus method (simply set $\alpha\defeq \vertices\of G\setminus\beta$).

%Instead of giving a direct proof for the % 130 years old
%identity \eqref{eq:tanner}, we state and prove
%generalizations.

%\subsection{Generalizations: Identities due to Wenzel and Krattenthaler}
%\label{sec:generalizations}
Tanner's identity (and the fact that it is of the \wpit) is obtained
immediately from the following generalization (and its bijective proof) which Hamel
\cite[Theorem~2.1]{hamel:2001} attributes to
%we present more general Pfaffian identi\-ties due
Ohta \cite{ohta:1993}; it was also found by Wenzel (see %also 
% \cite[Theorem~2.1]{hamel:2001},
\cite[Proposition 2.3]{wenzel:1993} and  \cite[Theorem 1]{dress-wenzel:1995}).

For arbitrary sets $A$ and $B$, denote the \EM{symmetric difference} of $A$ and $B$ by
$$
\symdef{A}{B}\defeq \pas{A\cup B}\setminus\pas{A\cap B}.
$$
%For fixed $n$, consider the complete
%graph $K_n$ with $\vertices\of{K_n}=\setof{1,\dots,n}$ and weight function $\weight$.
%We denote the Pfaffian corresponding
%to the subgraph of $K_n$ induced by some subset $Y\subseteq V\of{K_n}$ by
%$\pfaffian\of Y$ (see Definition~\ref{dfn:Stembridge}).

%The following generalization of Tanner's identity is equivalent to Equation~(5.0) in \cite{knuth:1996}, to Proposition 2.3 in \cite{wenzel:1993} and to Theorem~2.1 in
%\cite{hamel:2001} (according to Hamel \cite{hamel:2001} it is due to Ohta
%\cite{ohta:1993} and can also be found in \cite{hirota:1996}):
%We present it here in the ``Stembridge--style'':
\begin{thm}[Ohta]
\label{thm:wenzel}
Assume that the ordered set $\gamma=\pas{v_1,v_2,\dots,v_n}$ appears as $\gamma=\alpha\cup \beta$
% be two subsets of  % of odd cardinality
with
$\symdef{\alpha}{\beta}=\pas{v_{i_1},\dots,v_{i_t}}$.  Then we have:
\begin{equation}
\label{eq:wenzel}
\sum_{\tau=1}^t\pas{-1}^{\tau}\cdot
	\pfaffian\of{\symdef{\alpha}{\setof{v_{i_\tau}}}}\cdot
	\pfaffian\of{\symdef{\beta}{\setof{v_{i_\tau}}}} = 0.
\end{equation}
(Recall again that all subsets \EM{inherit} the order of $\word\gamma$.)

This identity is of the \wpit.
\end{thm}

%Before proving \eqref{eq:wenzel} we make two observations:
%\bit
%\item This identity directly implies
%Kuo's Proposition~\ref{prop:kuo}: Keep in mind the Kasteleyn--Percus method,
%set $\alpha=\vertices\of{G}\setminus\setof{a}$, $\beta=\vertices\of{G}\setminus\setof{b,c,d}$ and have a look at \figref{fig:kuo-by-wenzel}.
%\item This identity directly implies
%Tanner's identity~\eqref{eq:tanner}: Simply interpret the word $\alpha\, x$ as the
%ordered set $\alpha$, and the word $\pas{\alpha\, \beta}\setminus x$ as the ordered set
%$\beta$, and note that the signs in \eqref{eq:tanner} equal the signs in \eqref{eq:wenzel}
%in this interpretation: {\red Work this out!!!}
%\eit

% \psframebox[boxsep=false]{
%\begin{figure}
%\caption{Consider again Kuo's Proposition~\ref{prop:kuo} and observe the
%applicability of \eqref{eq:wenzel}. Set $\alpha=\vertices\of{G}\setminus\setof{a}$ and $\beta=\vertices\of{G}\setminus\setof{b,c,d}$
%in Ohta's Theorem~\ref{thm:wenzel}: In this special case, the four Pfaffians in 
%\eqref{eq:wenzel} correspond to the four terms in \eqref{eq:kuos-proposition};
%see also \figref{fig:matchings}.
%}
%\label{fig:kuo-by-wenzel}
%\input graphics/matchings-start
%\end{figure}
% }

The following proof given by Krattenthaler \cite{krattenthaler:2006} again uses superposition
of matchings, together with proper accounting for the sign--changes associated
to swapping colours in bicoloured paths. % , according to Observation~\ref{obs:arcsnsigns}.
We state this sign--change as follows:
\begin{lem}
\label{lem:swap-colours}
Consider the complete graph $K_V$ on the ordered set of vertices $V$, given as
disjoint union of red, blue and white vertices
$$V=\theblues\cup\thereds\cup\thewhites=\pas{v_1,v_2,\dots,v_n},$$
where
$\cardof{\theblues}\equiv\cardof{\thereds}\pmod 2$. Draw the edges of $K_V$ in
the specific way described
in Definition~\ref{dfn:Stembridge}. Let  $B$ be the corresponding bicoloured graph,
let $S={\blue\mu}\cup{\red\nu}$ be a superposition of matchings in $B$ (recall
Observation~\ref{obs:interpretation}), and let $p$ be a bicoloured path in $S$
(recall Observation~\ref{obs:nonintersecting}) with end vertices $v_i$ and $v_j$.
Swapping colours in $p$ (recall
Observation~\ref{obs:swap}) gives a superposition of matchings
$S={\blue\mu}^\prime\cup{\red\nu}^\prime$ in a bicoloured graph $B^\prime$, and we have
\begin{equation*}
\sgn\of{\blue\mu}\cdot\sgn\of{\red\nu} =
\pas{-1}^{\numof{\text{vertices of }\pas{\theblues\cup\thereds}\text{ between }v_i\text{ and }v_j}}
\sgn\of{\blue\mu^\prime}\cdot\sgn\of{\red\nu^\prime}.
\end{equation*}
If $v_i$ and $v_j$ appear in the ordered set of \EM{coloured}
vertices
$$
\thecoloureds=\theblues\cup\thereds=\pas{w_{1},w_{2},\dots,w_{{\cardof{\thecoloureds}}}}
$$
as $v_i=w_{{x}}$ and $v_{j}=w_{y}$, then this amounts to
\begin{equation}
\label{eq:sign-change-simple}
\sgn\of{\blue\mu}\cdot\sgn\of{\red\nu} =
\pas{-1}^{y-x+1}
\sgn\of{\blue\mu^\prime}\cdot\sgn\of{\red\nu^\prime}.
\end{equation}
\end{lem}
\begin{proof}
According to Observation~\ref{obs:arcsnsigns}, the change in sign corresponding to removing
some \EM{single} arc $\setof{v_k,v_l}$ from ${\blue\mu}$ and adding it to ${{\red\nu}}$ amounts
to
\begin{multline*}
\pas{-1}^{
	\numof{\text{vertices of }\pas{{\theblues\cup\thewhites}}\text{ between }v_k\text{ and }v_l} - 	\numof{\text{vertices of }\pas{{\thereds\cup\thewhites}}\text{ between }v_k\text{ and }v_l}
}
=\\
\pas{-1}^{\numof{\text{vertices of }\pas{\theblues\cup\thereds}\text{ between }v_k\text{ and }v_l}}.
\end{multline*}
Recolouring \EM{all} the arcs in path $p$ with end vertices $v_i$ and $v_j$ thus
gives a change in sign equal to the product
of \EM{all} these single sign--changes, which clearly amounts to
$$
\pas{-1}^{\numof{\text{vertices of }\pas{\theblues\cup\thereds}\text{ between }v_i\text{ and }v_j}}.
$$
\end{proof}
\begin{proof}[Proof of Theorem~\ref{thm:wenzel}]
For the combinatorial interpretation of the left--hand side of \eqref{eq:wenzel}, simply
combine Observation~\ref{obs:interpretation} (superposition of matchings) with % the
the definition of Pfaffians as given in Definition~\ref{dfn:Stembridge}: % (Pfaffians ``\`a la Stembridge''):
after expansion of the
products of Pfaffians, the typical summand is of the form
$$
\pas{-1}^\tau\cdot\sgn\of{{\blue\mu}}\cdot\sgn\of{{\red\nu}}\cdot\weight\of{{\blue\mu}}\cdot\weight\of{{\red\nu}},
$$
where $\pas{{\blue\mu},{\red\nu}}$ can be interpreted as superposition of matchings in the
bicoloured graph $\bicoloured{G}{\thereds}{\theblues}$ derived from the complete graph
$G$ on vertices $\alpha\cup \beta$, with
\bit
\item $\theblues=\symdef{\pas{\alpha\setminus \beta}}{\setof{v_{i_\tau}}}$,
\item $\thereds=\symdef{\pas{\beta\setminus \alpha}}{\setof{v_{i_\tau}}}$.
\eit
% Let $w=v_{i_\sigma}$.
According to Observation~\ref{obs:interpretation},
the operation $\swapcolours{v}$
of swapping colours in the unique path $p$ starting
in vertex $v=v_{i_\tau}$ is an \EM{involution} which preserves the weight
$\weight\of{{\blue\mu}}\cdot\weight\of{{\red\nu}}$.

So the proof is complete if we can show that $\swapcolours{v}$ is in fact
\EM{sign--reversing} in the following sense. Note that $\swapcolours{v}$ yields a superposition of
matchings $\pas{{\blue\mu}^\prime,{\red\nu}^\prime}$ in the
bicoloured graph $\bicoloured{G}{\thereds^\prime}{\theblues^\prime}$ with
\bit
\item $\theblues=\symdef{\pas{\alpha\setminus \beta}}{\setof{v_{i_\rho}}}$,
\item $\thereds=\symdef{\pas{\beta\setminus \alpha}}{\setof{v_{i_\rho}}}$,
\eit
where $v_{i_\rho}$ is the other endpoint of the unique path $p$. Thus we have to show
\begin{equation*}
%\label{eq:sign-reverse}
\pas{-1}^{\tau}\cdot\sgn\of{{\blue\mu}}\cdot\sgn\of{{\red\nu}}
=
-\pas{-1}^{\rho}\cdot\sgn\of{{\blue\mu^\prime}}\cdot\sgn\of{{\red\nu^\prime}}.
\end{equation*}
But this follows immediately from Lemma~\ref{lem:swap-colours}.
% Observation~\ref{obs:arcsnsigns}:
% (see also \figref{fig:arcsnsigns}): 
%the change in sign corresponding to removing
%some single arc $\setof{i,j}$ from ${\blue\mu}$ and adding it to ${{\red\nu}}$ amounts
%to
%$$
%\pas{-1}^{
%	\numof{\text{vertices of }\pas{{\alpha}\setminus{\beta}}\text{ between }i\text{ and }j} + 	\numof{\text{vertices of }\pas{{\beta}\setminus{\alpha}}\text{ between }i\text{ and }j}
%}.
%$$
%Recolouring \EM{all} the arcs in path $p$ thus gives a change in sign equal to the product
%of \EM{all} these single sign--changes. 
%Since we may assume  the
%vertices $\setof{i_1,\dots,i_t}$ to appear in ascending order, the change in
%sign simply amounts to
%$$
%\pas{-1}^{
%	\numof{\text{vertices of }\pas{\symdef{\alpha}{\beta}}\text{ between }\sigma\text{ and }\rho}
%}=\pas{-1}^{\sigma-\rho-1},
%$$
%which proves \eqref{eq:sign-reverse}.
\end{proof}

The same idea of proof applies to the following generalization, which to the best of my
knowledge is due to Krattenthaler \cite{krattenthaler:2006}:
\begin{thm}[Krattenthaler]
\label{thm:kratt1}
Assume that the ordered set $\gamma=\pas{v_1,v_2,\dots,v_n}$ appears as $\gamma=\alpha\cup \beta$. Let 
% be two subsets of  % of odd cardinality
% with
$M=\symdef{\alpha}{\beta}$.
%Let $\alpha$, $\beta$ be two subsets of $\pas{v_1,v_2,\dots,v_n}$ of the same cardinality $\pmod 2$
%(i.e., $\cardof{\alpha}\equiv\cardof{\beta}\pmod 2$), % of odd cardinality,
% with
% $\symdef{\alpha}{\beta}=\setof{i_1,\dots,i_t}$, $i_1<\cdots<i_t$.
% and let $M=\symdef{\alpha}{\beta}$.

If $\cardof{\alpha}$ is odd, then we have for all $s\geq 0$:
\begin{equation}
\label{eq:kratt1}
\sum_{%1\leq\tau_1<\cdots<\tau_{2s+1}\leq t}
	\substack{Y\subseteq M\\\cardof{Y}=2s+1}
	}
	\pas{-1}^{%\tau_1+\cdots+\tau_{2s+1}
		\SETSUM{Y}{M}
	}\cdot
	\pfaffian\of{\symdef{\alpha}{Y}}
		%\setof{i_{\tau_1},\dots,i_{\tau_{2s+1}}}}}
		\cdot
	\pfaffian\of{\symdef{\beta}{Y}}%\setof{i_{\tau_1},\dots,i_{\tau_{2s+1}}}}} \\ 
	= 0.
\end{equation}

If $\cardof{\alpha}$ is even, then we only have the weaker statement
\begin{equation}
\label{eq:kratt2}
\sum_{s=0}^{\floor{t/2}}
\sum_{%1\leq\tau_1<\cdots<\tau_{2s}\leq t}\pas{-1}^{\tau_1+\cdots+\tau_{2s}
	\substack{Y\subseteq M\\\cardof{Y}=2s}
	}
	\pas{-1}^{
		\SETSUM{Y}{M}
	}\cdot
	\pfaffian\of{\symdef{\alpha}{Y}}\cdot
	\pfaffian\of{\symdef{\beta}{Y}} = 0.
\end{equation}

These identities are of the \wpit.
\end{thm}
\begin{proof}
For every superposition of matchings $\pas{{\blue\mu},{\red\nu}}$ involved
in \eqref{eq:kratt1} (in the sense of the proof of Theorem~\ref{thm:wenzel}) consider
the subset $S\subseteq Y$ of vertices $v$
with the property that the \EM{other} endpoint of the unique path starting in $v$
\EM{does not} belong to $Y$. Note that $S$ is of
\EM{odd} cardinality, so in particular $S\neq\emptyset$.

By the same reasoning as in the proof of Theorem~\ref{thm:wenzel}, \EM{simultaneously}
swapping colours in \EM{all} paths starting in vertices from $S$ yields a sign--reversing
and weight--preserving involution, which proves \eqref{eq:kratt1}.

Now consider the family of superpositions of matchings
corresponding to \eqref{eq:kratt2}. Let $v$ be some arbitrary, but fixed element in $M$,
and note that the operation
$\swapcolours{v}$ of swapping colours in the unique
path with starting point $v$ yields % a sign--reversing and
weight--preserving involution, which is sign--reversing according to Lemma~\ref{lem:swap-colours}.
\end{proof}

Since $\pfaffian\of{X}\equiv0$ if the cardinality of $X$ is odd, we may restate
(a weaker version of) the above theorem in a uniform way. 
% for some arbitrary subset of $X\subset\N$.
\begin{cor}
\label{thm:kratt3}
Assume that the ordered set $\gamma=\pas{v_1,v_2,\dots,v_n}$ appears as $\gamma=\alpha\cup \beta$. Let 
$M=\symdef{\alpha}{\beta}$.
%Let $\alpha$, $\beta$ be two subsets of $\setof{v_1,v_2,\dots,v_n}$, let $M=\symdef{\alpha}{\beta}$.
 % of the same cardinality $\pmod 2$.
% with $\symdef{\alpha}{\beta}=\setof{i_1,\dots,i_t}$, $i_1<\cdots<i_t$.
Then we have:

\begin{equation}
\label{eq:kratt3}
\sum_{Y\subseteq M}
\pas{-1}^{\SETSUM{Y}{M}}\cdot
	\pfaffian\of{\symdef{\alpha}{Y}}\cdot
	\pfaffian\of{\symdef{\beta}{Y}} = 0.
\end{equation}
This identity is of the \wpit.
\end{cor}

\subsection{Further applications}
% The same ideas yield the following results.
%For the rest of this paper,
%consider the complete graph $G=K_V$ with (ordered) vertex set $V=\pas{v_1,v_2,\dots,v_n}$,
%where $V$ appears as the disjoint union of blue, red and white vertices:
%$V=\markextra{\theblues}\cup\markextra{\thereds}\cup\thewhites$.
%Denote the set of coloured vertices by
%$\thecoloureds=\markextra{\theblues}\cup\markextra{\thereds}$.
For the rest of this paper,
consider the complete graph $G=K_V$ with (ordered) vertex set $V=\pas{v_1,v_2,\dots,v_n}$,
where $V$ appears as the disjoint union of coloured vertices and white vertices,
$V=\markextra{\thecoloureds}\disjuni\thewhites$, and where the set of coloured vertices
% is of even cardinality and 
is partitioned into two sets of equal size $\thecoloureds=\thelefts\disjuni\therights$,
$\cardof\thelefts=\cardof\therights$. For the following assertions, imagine that
\EM{all} vertices in $\therights$
are ``initially'' blue and that \EM{all} vertices in $\thelefts$ are ``initially'' red,
and bicoloured graphs are derived therefrom by changing this ``initial'' colouring.

\subsubsection{Srinivasan's result}

For the next Lemma, let $X$ be some fixed subset
$X\subseteq\therights$, and let $\theblues\defeq\therights\setminus X$ and
$\thereds\defeq\thelefts\cup X$. Consider the set of
all superpositions of matchings ${{\red\mu}\cup{\blue\nu}}$
in $\bicoloured{{G}}{{{\markextra{\thereds}}}}{{{\markextra{\theblues}}}}$
with the \EM{additional} property,
that every bicoloured path in ${{\red\mu}\cup{\blue\nu}}$ has \EM{at most one}
end vertex in $\markextra{\therights}$: denote this set by $\strangeset$.
% {\markextra{\theblues}}{X}$.
To every object % $\pas{{\red\mu},{\blue\nu}}\in{\mathcal F}_X$ 
${{\red\mu}\cup{\blue\nu}}\in\strangeset$
%{\markextra{\theblues}}{X}$
assign signed weight
$$
\weight\of{{\red\mu}\cup{\blue\nu}}\defeq
\pas{-1}^{\SETSUM{B\setminus X}{\thecoloureds}}
%\pas{-1}^{\SETSUM{X}{\thecoloureds}}
\sgn\of{{\red\mu}}\cdot\sgn\of{{\blue\nu}}\cdot
\weight\of{{\red\mu}}\cdot\weight\of{{\blue\nu}}.
$$
Now
% assume two subsets $A$, $B$ (not necessarily disjoint) of $V$ with $V=A\cup B$.
% Set $\thecoloureds\defeq\symdef{A}{B}$, $\markextra{\thereds}=A\setminus B$ and
% $\markextra{\theblues}=B\setminus A$,
% and
consider the generating function $F$ of
$\bigcup_{X\subseteq\markextra{\therights}}{\strangeset}$,
i.e.,
$$
F\defeq\sum_{X\subseteq\markextra{\therights}}\;
	\sum_{f\in\strangeset}\weight\of f.$$

\begin{lem}
\label{lem:srinivasan}
With the above definitions, we have for the generating function $F$:
\begin{equation}
\label{eq:general-srinivasan}
F = 
\pas{-1}^{\SETSUM{B}{\thecoloureds}}\cdot
\sum_{X\subseteq\markextra{\therights}}
	\pas{-1}^{\SETSUM{X}{\thecoloureds}}\cdot
	\pfaffian\of{\markextra{\thelefts}\cup\thewhites\cup X}\cdot
	\pfaffian\of{\pas{\markextra{\therights}\cup\thewhites}\setminus X}.
\end{equation}
\end{lem}
\begin{proof}
Clearly, all the signed weights of objects from
$\bigcup_{X\subseteq\markextra{\therights}}{\strangeset}$
\EM{do appear}
in the sum on the right--hand side of \eqref{eq:general-srinivasan}.
So we have to show that all the ``superfluous'' terms in \eqref{eq:general-srinivasan}
cancel.

These ``superfluous'' objects are precisely superpositions of matchings, where there
\EM{exists} a bicoloured path connecting \EM{two} vertices of $\markextra{\therights}$:
of all such paths 
choose the one, $p$, with the smallest end vertex, and swap colours in $p$. This gives
a weight--preserving involution, which is sign--reversing according to Lemma~\ref{lem:swap-colours}.
\end{proof}
From this, we easily deduce (a slight generalization of) a result of Srinivasan
(\cite[Corollary 3.2]{srinivasan:1994}, see also \cite[Theorem 3.3]{hamel:2001}):

\begin{cor}[Srinivasan]
Let $m+n$ be an even integer and consider the ordered set $V=\pas{v_1,\dots,v_{m+n}}$, partitioned into the disjoint subsets $A=\pas{v_{i_1},\dots,v_{i_m}}$ and $B=\pas{v_{j_1},\dots,v_{j_n}}$.
%Let $A$, $B$ be two \EM{disjoint} sets, $\cardof{A}=m$, $\cardof{B}=n$ with
%$\pas{m+n}$ even. Denote $M=A\cup B$. % of even (???the same???) cardinality.
%% with $\cardof{M}$ even. % (???)
Then we have
the following expansions:

If $m<n$, then
\begin{equation}
\label{eq:srinivasan_p_less_q}
\pfaffian\of{V} = -\sum_{\substack{X\subset B,\\X\neq B}}
\pas{-1}^{\SETSUM{B\setminus X}{V}}\cdot
	\pfaffian\of{A\cup X}\cdot
	\pfaffian\of{B\setminus X}.
\end{equation}

If $m=n$, then
\begin{equation}
\label{eq:srinivasan_p_equals_q}
\pfaffian\of{V}= \pfaffian\of{A,B}-\sum_{\substack{X\subset B,\\X\neq B}}
\pas{-1}^{\SETSUM{B\setminus X}{V}}\cdot
	\pfaffian\of{A\cup X}\cdot
	\pfaffian\of{B\setminus X}.
\end{equation}

If $m>n$, then
\begin{multline}
\label{eq:srinivasan_p_greater_q}
\pfaffian\of{V}=
\pas{-1}^{\binom{n+1}{2}}
\sum_{\substack{Y\subset A,\\\cardof{Y}=n}}
	\pas{-1}^{\SETSUM{Y}{A}}
	\pfaffian\of{A\setminus Y}\cdot\pfaffian\of{B, Y}
\\ -
\sum_{\substack{X\subset B,\\X\neq B}}
\pas{-1}^{\SETSUM{B\setminus X}{V}}\cdot
	\pfaffian\of{A\cup X}\cdot
	\pfaffian\of{B\setminus X}.
\end{multline}
\end{cor}
\begin{proof}
We apply Lemma~\ref{lem:srinivasan} with $\markextra{\thelefts}=A$, 
$\markextra{\therights}=B$ and $\thewhites=\emptyset$: since there are no white vertices,
% we have
%$\markextra{\thereds}=A$, $\markextra{\theblues}=B$ and
%$\thecoloureds=M$.
%% and $\alpha\cap \beta=\emptyset$.
%Note that according to \eqref{eq:general-srinivasan}, the sum substracted in the % over the subsets $X$, which appears in the
%\RHS\ of
%equations~\eqref{eq:srinivasan_p_less_q},
%\eqref{eq:srinivasan_p_equals_q} and
%\eqref{eq:srinivasan_p_greater_q}
%equals $F-\pfaffian\of{A\cup B}$.
%Moreover,
bicoloured paths in \EM{any} superposition of matchings
$\pas{{\red\mu}\cup{\blue\nu}}\in\bigcup_{X\subseteq B}{\strangeset}$
are simply \EM{edges},
whose end vertices are of the \EM{same} colour: hence % for an object from $\mathcal F$,
we must have $X=B$ (because there is no matching
on $B\setminus X\neq\emptyset$ where every edge has only \EM{one} end vertex in $B$).
So ${{\red\mu}\cup{\blue\nu}}$ corresponds to a \EM{unique} 
matching in the semi--bipartite graph $\semibip{B}{A}$, and vice versa, whence
we have
$$
F = \pfaffian\of{\semibip{B}{A}}.
$$
Expanding the trivial identity
$$
\pfaffian\of{V} = \pfaffian\of{\semibip{B}{A}}
- \pas{F- \pfaffian\of{A\cup B}\cdot\pfaffian\of\emptyset}
$$
according to \eqref{eq:proof-semibip} and \eqref{eq:general-srinivasan}, respectively,
proves all three assertions: note that $\pfaffian\of{\semibip{B}{A}}=0$ if $m<n$, and
$\pfaffian\of{\semibip{B}{A}}=\pfaffian\of{A,B}$ if $m=n$.
\end{proof}

\subsubsection{Graphical condensation}
For the rest of this paper, assume that ${\therights}=\pas{b_1,b_2\dots,b_k}$
and ${\thelefts}=\pas{r_1,r_2,\dots,r_k}$ and that the ordered set of
coloured vertices appears as $\thecoloureds=\pas{r_1,b_1,r_2,b_2\dots,r_k,b_k}$.
\begin{dfn}[Planar weight function]
We call a weight function $\weight$ on $\edges\of{G}$ a \EM{planar weight function}
if it assigns weight \EM{zero} to % has the property that 
\EM{every} superposition of matchings ${{\red\mu}\cup{\blue\nu}}$ in
\EM{every} $\bicoloured{{G}}{\thereds}{\theblues}$ (where
$\thereds\disjuni\theblues=\thelefts\disjuni\therights$ is an \EM{arbitrary} partition
of the coloured vertices in two subsets of blue and red vertices)
 that contains a
bicoloured path $p$ connecting two vertices from $\therights$ or two vertices from $\thelefts$.
%has weight \EM{zero}:
%$\weight\of{{\red\mu}\cup{\blue\nu}}= 0$.
\end{dfn}
\begin{lem}
\label{lem:sign-preserving}
In the setting described above, assume that the weight function is \EM{planar}.
% definition of $G=K_V$, $V={\therights\cup\thelefts\cup\thewhites}=\pas{b_1,b_2\dots,b_k}$
%and ${\thereds}=\pas{r_1,r_2,\dots,r_k}$.
%Furthermore, assume
%\bit
%\item that the ordered set of coloured vertices appears as
%	$\thecoloureds=\pas{r_1,b_1,r_2,b_2\dots,r_k,b_k}$,
%%	where the ordered sets of red and blue vertices appear as
%%	$\thereds=\pas{r_1,r_2,\dots,r_k}$
%%	and
%%	$\theblues=\pas{b_1,b_2\dots,b_k}$,
%%	respectively,
%\item and that the weight function $\weight$ on
%$\edges\of{K_V}$ has the property that 
%%\begin{quote}
%%For all $v,w\in\thereds$, every bicoloured path $p$ connecting $v$ and $w$ in (the superpositions of
%%matchings corresponding to) $\pfaffian\of{\alpha}\cdot \pfaffian\of{\beta}$ has weight zero:
%%$\weight\of p = 0$.
%\EM{every} superposition of matchings $\pas{{\red\mu}^\prime\cup{\blue\nu}^\prime}$ in
%\EM{every} $\bicoloured{{G}}{\thereds^\prime}{\theblues^\prime}$ (where
%$\thereds^\prime\disjuni\theblues^\prime=\thereds\disjuni\theblues$)
% that contains a
%bicoloured path $p$ connecting two vertices from $\theblues$ or two vertices from $\thereds$
%% whose end vertices are of the \EM{same} colour (i.e., both red or both blue)
%has weight \EM{zero}:
%$\weight\of{{\red\mu}^\prime\cup{\blue\nu}^\prime}= 0$.
%\eit
%%\end{quote}

Then we have for all (fixed) subsets $X\subseteq\thelefts$:
\begin{multline}
\label{eq:sign-preserving}
\sum_{W\subseteq\therights}
	\pfaffian\of{\pas{\thelefts\cup\thewhites}\cup W}
	\cdot
	\pfaffian\of{\pas{\therights\cup\thewhites}\setminus W} \\
=
\sum_{V\subseteq\therights}
	\pfaffian\of{\pas{\pas{\thelefts\cup\thewhites}\cup V} \setminus X}
	\cdot
	\pfaffian\of{\pas{\pas{\therights\cup\thewhites}\setminus V}\cup X}.
\end{multline}
This identity is of the \wpit.
\end{lem}
\begin{proof}
The left--hand side of
\eqref{eq:sign-preserving}
corresponds to superpositions of matchings ${{\red\mu}\cup{\blue\nu}}$
\EM{of non--vanishing weight} % ($\weight\of{{\red\mu}}\cdot\weight\of{{\blue\nu}}\neq0$),
with red vertices $\thelefts\cup W$
and blue vertices $\therights\setminus W$, for some
subset $W\subseteq\therights$.
The assumption implies that every bicoloured path $p$,
%of non--vanishing weight ($\weight\of p\neq 0$),
which starts in some vertex $v\in\thelefts$, must have
its other end in $\therights$. So %Now note that
according to %\eqref{eq:sign-reverse} 
\eqref{eq:sign-change-simple}
(and the
specific ``alternating ordering'' of $\thecoloureds=\therights\cup\thelefts$), the operation $\swapcolours{v}$ of swapping
colours in $p$ gives
a weight-- \EM{and} sign--preserving bijection. Thus, swapping colours in \EM{all}
bicoloured paths
with end vertex in the fixed subset
$X=\setof{v_1,\dots,v_m}\subseteq\thelefts$ % (note that this end vertex is red!)
gives a weight-- and sign--preserving
bijection $\chi\defeq\swapcolours{v_1}\concat\dots\concat\swapcolours{v_m}$.

Denote by $Y$ the set of the \EM{other} end vertices of these paths. Note that
$Y\subseteq\therights$ by assumption and consider the set $V\defeq\symdef{W}{Y}$: it is obvious that the image of $\chi$
% whose image
is precisely the set of superpositions of matchings corresponding to the
right hand side of \eqref{eq:sign-preserving}.
\end{proof}

This yields immediately the following generalization of a result by Yan, Yeh and Zhang
\cite[Theorem~2.2]{yan-yeh-zhang:2005} (also given in Kuo \cite[Theorem 2.1 and Theorem 2.3]{kuo:2006}):
\begin{cor}
\label{cor:kuos-theorem}
Let $G$ be a \EM{planar} graph with vertices $r_1,b_1,\dots,r_k,b_k$ 
appearing in that cyclic order on the boundary of a face of $G$. Let $\thelefts=\setof{r_1,\dots,r_k}$ and $\therights=\setof{b_1,\dots,b_k}$. Then we have for every fixed subset $X\subseteq\thelefts$
\begin{equation*}
\sum_{W\subseteq \theblues}
	\nofmatchings{\subgraph{G}{\pas{\therights\setminus W}}}
	\nofmatchings{\subgraph{G}{\pas{\thelefts \cup {W}}}}
=
\sum_{V\subseteq \therights}
	\nofmatchings{\subgraph{G}{\pas{\pas{\therights\setminus V} \cup X}}}
	\nofmatchings{\subgraph{G}{\pas{\pas{\thelefts\cup V} \setminus X}}}.
\end{equation*}
\end{cor}
\begin{proof}
Recall that no two different paths arising from some superposition of matchings can have a
vertex in common. So the ordering of $\thelefts\cup\therights$ along the boundary of a face
of the \EM{planar} graph $G$ implies that no bicoloured path can have \EM{both} end vertices
in $\thelefts$ or in $\therights$, hence the weight function is \EM{planar} and we may
apply Lemma~\ref{lem:sign-preserving}. By the Kasteleyn--Percus method, the assertion follows
from \eqref{eq:sign-preserving}.
%(otherwise, it would intersect some other bicoloured path
\end{proof}

\subsubsection{Ciucu's matching factorization theorem}
In the same manner, we may prove Ciucu's matching factorization Theorem
\cite[Theorem~1.2]{ciucu:1997},
in the (equivalent) formulation given in \cite[Theorem~2.2]{yan-zhang:2005}.

\begin{lem}
\label{lem:ciucu}
In the setting described above, assume that the weight function is \EM{planar}. Moreover,
assume that the set of coloured vertices % $\thecoloureds=\theleft\disjuni\therights$
is partitioned into two subsets of equal size, $\thecoloureds=U\disjuni V$,
$\cardof U = \cardof V$, such 
\bit
\item that there is
\EM{no even--length path} of non--vanishing weight connecting \EM{any two}
vertices $\pas{x,y}$ with $x\in U$ and $y\in V$,
\item and that there is
\EM{no odd--length path} of non--vanishing weight connecting \EM{any two}
vertices $\pas{x,y}$ with $x, y\in U$ or $x, y\in V$.
\eit
Let $\therights_U\defeq\therights\cap U$, $\therights_V\defeq\therights\cap V$, 
$\thelefts_U\defeq\thelefts\cap U$ and $\thelefts_V\defeq\thelefts\cap U$.
(Note that $\cardof{\therights_U}=\cardof{\thelefts_V}$ and
$\cardof{\therights_V}=\cardof{\thelefts_U}$.)
Then we have
\begin{multline}
\label{eq:ciucu}
2^k\cdot
\pfaffian\of{\thewhites\cup \therights_U\cup\thelefts_V}\cdot
\pfaffian\of{\thewhites\cup \therights_V\cup\thelefts_U}
=\\
\sum_{\substack{X\subseteq V,Y\subseteq U\\\cardof{X}=\cardof{Y}}}
\pfaffian\of{\thewhites\cup X \cup Y}\cdot
\pfaffian\of{\thewhites\cup \complement{X}\cup\complement{Y}}.
\end{multline}
This identity is of the \wpit.
\end{lem}
\begin{proof}
First, consider an \EM{arbitrary} superposition of matchings % ${\red\mu}\cup{\blue\nu}$
in the bicoloured graph $\bicoloured{G}{\theblues}{\thereds}$, where
with $\theblues\disjuni\thereds$ is an \EM{arbitrary} partition of the set
$\thelefts\disjuni\therights$. Since
\bit
\item bicoloured paths of \EM{even} length have end vertices
of \EM{different} colour,
\item and bicoloured paths of \EM{odd} length have end vertices of the
\EM{same} colour,
\eit
the assumption implies that if $\weight\of{{\red\mu}\cup{\blue\nu}}\neq 0$, then
every bicoloured path in ${\red\mu}\cup{\blue\nu}$ must connect (vertices from)
$\thereds\cap U$ with $\theblues\cap U$, $\theblues\cap U$ with $\theblues\cap V$, 
$\theblues\cap V$ with $\thereds\cap V$, or $\thereds\cap V$ with $\thereds\cap U$.
See the left
picture in \figref{fig:small_table}, where the possible connections by bicoloured paths are
indicated by arrows.

Now consider the product of Pfaffians in the left hand side of \eqref{eq:ciucu}, i.e., choose $\theblues = \therights_U\cup\thelefts_V$ and
$\thereds  = \therights_V\cup\thelefts_U$. Note that in this case we have
$\thereds\cap U=\thelefts_U$, $\thereds\cap V=\therights_V$,
$\theblues\cap U=\therights_U$ and $\theblues\cap V=\thelefts_V$, and the
assertion about the possible connections with bicoloured paths (as depicted
in the left picture of \figref{fig:small_table}) would hold also \EM{without} the
additional assumption, since the weight function $\weight$ is \EM{planar}.
 
%Since bicoloured paths of \EM{even} length have end vertices
%of \EM{different} colour, and bicoloured paths of \EM{odd} length have end vertices of the
%\EM{same} colour,
%the assumption implies that if $\weight\of{{\red\mu}\cup{\blue\nu}}\neq 0$, then
%every bicoloured path in ${\red\mu}\cup{\blue\nu}$ must connect (vertices from)
%$\thelefts_U$ with $\therights_U$, $\therights_U$ with $\thelefts_V$, 
%$\thelefts_V$ with $\therights_V$, or $\therights_V$ with $\thelefts_U$ (see the left
%picture in \figref{fig:small_table}, where the possible connections by bicoloured paths are
%indicated by arrows).

\begin{figure}
\caption{Illustration to the proof of Lemma~\ref{lem:ciucu}: The left picture shows
the partition of the set of coloured vertices
in four subsets, induced by the bipartitions  $\thereds\disjuni\theblues=U\disjuni V$.
The arrows indicate the \EM{only possible} connections by bicoloured paths: for instance,
no bicoloured path can connect a vertex from $\thereds\cap U$ with $\theblues\cap V$, since
such paths must have \EM{even} length, which is ruled out by the additional assumption in
Lemma~\ref{lem:ciucu}. The right picture shows the partition of the set of coloured
vertices in eight subsets, induced by the bipartitions
$\thereds\disjuni\theblues=U\disjuni V=\thelefts\disjuni\therights$.
The arrows indicate the \EM{only possible} connections between
``wrongly coloured'' vertices (these are the ones indicated by the inscribed square
with the dashed boundary) by bicoloured paths: in particular, every bicoloured path must
have either no ``wrongly coloured'' end vertex or two ``wrongly coloured'' end vertices.}
\label{fig:small_table}
\begin{center}
\input graphics/small_table2.tex
\end{center}
\end{figure}

Consider some arbitrary superposition of matchings ${\red\mu}\cup{\blue\nu}$
in $\bicoloured{G}{\theblues}{\thereds}$. Swapping colours in an arbitrary
subset of bicoloured paths in ${\red\mu}\cup{\blue\nu}$ gives a unique superposition of matchings ${\red\mu}^\prime\cup{\blue\nu}^\prime$ of the same
weight and sign, which appears in the right hand side of \eqref{eq:ciucu}. So let the
bicoloured paths $\pas{p_1,p_2,\dots,p_k}$ of ${\red\mu}\cup{\blue\nu}$
% in any superposition of matchings
appear in the order implied by their smaller
end vertex (in the order of the set of coloured vertices $\thecoloureds$)
and consider the set of pairs $\pas{S,{\red\mu}\cup{\blue\nu}}$, where $S=\setof{i_1,\dots,i_m}$ is a subset
of $\myrange{k}$. Clearly, this set of pairs corresponds to the left hand side of
\eqref{eq:ciucu},
% and ${\red\mu}\cup{\blue\nu}$ is a superposition of matchings 
% in the bicoloured graph $\bicoloured{G}{\theblues}{\thereds}$:
and swapping colours in the bicoloured paths $p_{i_1},\dots,p_{i_m}$ defines a weight--
and sign--preserving injective mapping $\chi$ into the set of \soms\ corresponding to
the right hand side of \eqref{eq:ciucu}.

So it remains to show that this mapping is in fact surjective.
To construct the preimage $\chi^{-1}\of{{\red\mu}^\prime\cup{\blue\nu}^\prime}$ of some arbitrary superposition of matchings
%$\pas{S,{\red\mu}^\prime\cup{\blue\nu}^\prime}$
corresponding to the right hand side of \eqref{eq:ciucu}, we would like to reverse the
operation $\chi$, i.e.,
\bit
\item note down the subset $S$ of paths (in the order described above) with ``wrongly
	coloured'' end vertices (these
	are vertices in $\pas{\therights\cap U} \cup\pas{\thelefts\cap V}$ which are red
	and vertices in $\pas{\therights\cap V} \cup\pas{\thelefts\cap U}$ which are blue),
\item and swap colours in all paths from $S$,
	\EM{hoping} that this operation gives some superposition of matchings ${\red\mu}\cup{\blue\nu}$
	corresponding to the left hand side of \eqref{eq:ciucu}.
\eit
This simple idea clearly will work if there is \EM{no} path with \EM{only one}
``wrongly coloured'' end vertex. To see that this is the case, just have a look at the right picture
in \figref{fig:small_table}: Since the weight function is \EM{planar}, it is not
possible for a superposition of matchings (of non--vanishing weight) to have
\EM{only one} ``wrongly coloured'' end vertex (for instance, there is no bicoloured connection
from the ``wrongly coloured'' set $\complement{X}\cap\therights$ to $\complement{Y}\cap\therights$ or to
$X\cap\therights$). Thus $\chi$ is in fact a weight--preserving and sign--preserving
bijection.
\end{proof}

\begin{cor}
\label{cor:ciucu}
Let $G$ be a \EM{planar} \EM{bipartite} graph, where the bipartition of the
vertex set of $G$ is given as $\vertices\of{G}=A\disjuni B$. Assume vertices $r_1,b_1,\dots,r_k,b_k$ 
appear in that cyclic order on the boundary of a face of $G$, and let $\thelefts=\setof{r_1,\dots,r_k}$ and $\therights=\setof{b_1,\dots,b_k}$. Let
$U\defeq A\cap\pas{\thelefts\disjuni\therights}$ and
$V\defeq B\cap\pas{\thelefts\disjuni\therights}$, and
assume that 
$\cardof{U} = \cardof{V}$. 

Let $\therights_U\defeq\therights\cap U$, $\therights_V\defeq\therights\cap V$, 
$\thelefts_U\defeq\thelefts\cap U$ and $\thelefts_V\defeq\thelefts\cap U$.
Then we have
\begin{equation*}
2^k\cdot
\nofmatchings{\subgraph{G}{\pas{\therights_U\cup\thelefts_V}}}\cdot
\nofmatchings{\subgraph{G}{\pas{\therights_V\cup\thelefts_U}}}
=\\
\sum_{\substack{X\subseteq V,Y\subseteq U\\\cardof{X}=\cardof{Y}}}
\nofmatchings{\subgraph{G}{\pas{X \cup Y}}}\cdot
\nofmatchings{\subgraph{G}{\pas{\complement{X}\cup\complement{Y}}}}.
\end{equation*}
\end{cor}
\begin{proof}
Note that the assumption that $G$ is \EM{bipartite} implies the additional
assumption of Lemma~\ref{lem:ciucu}.
The assertion now follows from \eqref{eq:ciucu} by the Kasteleyn--Percus method.
\end{proof}
\subsubsection{Graphical edge--condensation}

\begin{figure}
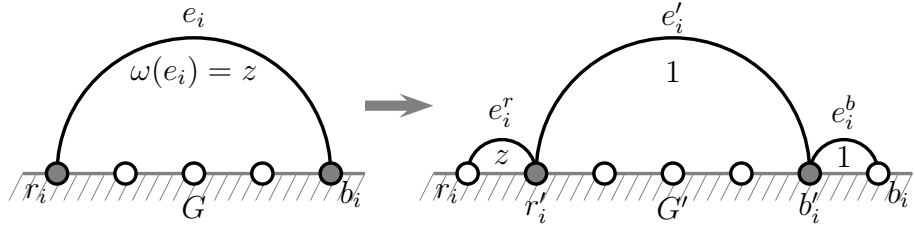

\caption{Illustration: $G^\prime$ is obtained from $G$ by ``subdividing'' edge $e_i$,
i.e., by replacing edge $e_i$ by the path $\pas{r_i, r_i^\prime, b_i^\prime, b_i}$,
thus introducing new vertices $r_i^\prime$ and $b_i^\prime$, and new edges
$e_i^r\defeq\setof{r_i, r_i^\prime}$, $e_i^\prime\defeq\setof{r_i^\prime, b_i^\prime}$ and $e_i^b\defeq\setof{b_i^\prime, b_i}$.}
\label{fig:edge-condensation}
\input graphics/edge-condensation
\end{figure}

Consider again the situation of Lemma~\ref{lem:sign-preserving}.
% and assume 
%that in the set $\vertices$ of ordered vertices, $r_i$ is immediately followed by $b_i$
%for all $i=1,\dots,k$ (see the left picture in \figref{fig:edge-condensation}).
Denote the edge
connecting vertex $r_i$ with vertex $b_i$ by $e_i\defeq\setof{r_i,b_i}$.
From the given complete
graph $G=K_V$ construct a new graph $G^\prime$ by replacing every edge $e_i$ by a path of length three $\pas{r_i,r_i^\prime,b_i^\prime,b_i}$, thus
inserting
\bit
\item new vertices $r_i^\prime$ and $b_i^\prime$, where $r_i^\prime$ is the
\EM{immediate successor} of $r_i$ and $b_i^\prime$ is the \EM{immediate predecessor}
of $b_i$, in the set of ordered vertices of $G^\prime$,
\item and new edges
$e_i^r\defeq\setof{r_i,r_i^\prime}$,
$e_i^\prime\defeq\setof{r_i^\prime,b_i^\prime}$
and $e_i^b\defeq\setof{b_i^\prime,b_i}$, where the weights of the new edges are given as
$\weight\of{e_i^r} = \weight\of{e_i}$ and
$\weight\of{e_i^\prime} =
\weight\of{e_i^b}=1$.
\eit
Note that $G^\prime$ is not a complete graph ($r_i^\prime$ and $b_i^\prime$ are vertices of degree 2 in $G$), but we may view it as a complete graph by introducing edges of weight zero.
(\figref{fig:edge-condensation} illustrates this construction.)

Observe
%\bit
% \item
that there is an obvious  \EM{weight--preserving} bijection between matchings of $G$
	and matchings of $G^\prime$: For aribtrary $\mu\in\matchings{G}$ and all $i=1,\dots, k$,
	% consider the edges $e_i$ and
	\bit
	\item \EM{replace} $e_i$ by edges $e_i^r$ and $e_i^b$ if $e_i\in\mu$,
	\item \EM{add} $e_i^\prime$ to $\mu$ if $e_i\not\in\mu$,
	\eit
	to obtain a matching $\mu^\prime\in\matchings{G^\prime}$.
% \item
Note that this bijection yields a change in sign equal to
$$
\prod_{i=1}^k\pas{-1}^{\numof{\text{vertices between $r_i$ and $b_i$}}} =
\pas{-1}^{k + \SETSUM{\thecoloureds}{V}},
$$
according to Observation~\ref{obs:arcsnsigns} (since
	the edges $e_i^r$ and $e_i^b$ in $G^\prime$
	can never be involved in any \EM{crossing}).
% \eit
Thus we obtain immediately
\begin{equation}
\label{eq:ciucu1996}
\pas{-1}^{k + \SETSUM{\thecoloureds}{V}}\pfaffian\of{G} = \pfaffian\of{G^\prime}.
\end{equation}
%In the same way we observe the following: C
Now consider the bicoloured graph $B^\prime=\bicoloured{G^\prime}{\theblues^\prime}{\thereds^\prime}$,
where $\theblues^\prime\disjuni\thereds^\prime$ is some partition of the set of coloured
vertices $\thecoloureds^\prime\defeq\setof{r_1^\prime,b_1^\prime,\dots, r_k^\prime,b_k^\prime}$
in $G^\prime$.
%, and observe:

%\theblues^\prime\disjuni\thereds^\prime=\thhcolureds = any bicoloured graph $B$ derived
%from $G^\prime$ with
%coloured vertices $\setof{r_1^\prime,b_1^\prime,\dots, r_k^\prime,b_k^\prime}$.
%\bit
%\item
If $r_i^\prime$ and $b_i^\prime$ are of \EM{different} colours, then edges $e_i^r$ and
$e_i^b$ must \EM{both} belong to \EM{every} superposition of matchings
${\red\mu}^\prime\cup{\blue\nu}^\prime$ in
the bicoloured graph $B^\prime=\bicoloured{G^\prime}{\theblues^\prime}{\thereds^\prime}$.
% So after accounting for 
%Since their total contribution
%$\weight\of{e_i^b}\cdot\weight\of{e_i^r}$ % \weight\of{e_i}$$
%to the weight of ${\red\mu}^\prime\cup{\blue\nu}^\prime$ equals $1$,
So we may simply remove them
together with their end--vertices from the corresponding subgraph (blue or red, respectively)
of $G^\prime$,
% $r_i^\prime$ and $b_i^\prime$,
thus obtaining a superposition of
matchings ${\red\mu}^{\prime\prime}\cup{\blue\nu}^{\prime\prime}$ 
%with %where
% with 
%in the remaining graph
%$\subgraph{G^\prime}{\setof{r_i^\prime,b_i^\prime}}$, with
%, where the ``surviving'' endpoints of the
%removed edges are of the opposite colour (i.e., $r_i$ is red if $e_i^r$ was blue, and vice
%versa, see \figref{fig:edge-condensation2} for
%an illustration). Altogether this implies
%$
%\weight\of{{\red\mu}^{\prime}\cup{\blue\nu}^{\prime}}
%=\weight\of{e_i}\times
%\weight\of{{\red\mu}^{\prime\prime}\cup{\blue\nu}^{\prime\prime}}
%$
%where ${\red\mu}^{\prime\prime}\cup{\blue\nu}^{\prime\prime}$ is a superposition of matchings
in the bicoloured graph
$
\bicoloured{\subgraph{G^\prime}{\setof{r_i^\prime,b_i^\prime}}}%
{\theblues^{\prime\prime}}{\thereds^{\prime\prime}}
$
with
\bit
\item 
$\theblues^{\prime\prime}\defeq\theblues^\prime\cup\setof{b_i}\setminus\setof{r_i^\prime}$
and 
$\thereds^{\prime\prime}\defeq\thereds^\prime\cup\setof{r_i}\setminus\setof{b_i^\prime}$,
if $r_i^\prime\in\theblues^\prime$ (as in \figref{fig:edge-condensation2}),
\item
$\theblues^{\prime\prime}\defeq\theblues^\prime\cup\setof{r_i}\setminus\setof{b_i^\prime}$
and 
$\thereds^{\prime\prime}\defeq\thereds^\prime\cup\setof{b_i}\setminus\setof{r_i^\prime}$,
if $r_i^\prime\in\thereds^\prime$.
\eit
Observe that
$\subgraph{G^\prime}{\setof{r_i^\prime,b_i^\prime}}$
``locally looks like'' the original graph $G$. For later use, note that
\begin{equation}
\label{eq:weight-change}
\weight\of{{\red\mu}^{\prime}\cup{\blue\nu}^{\prime}}
=\weight\of{e_i}\times
\weight\of{{\red\mu}^{\prime\prime}\cup{\blue\nu}^{\prime\prime}}
\end{equation}
(there is no sign--change here, since the edges $e_i^r$ and $e_i^b$ in $G^\prime$
can never be involved in any \EM{crossing}).
\begin{figure}
\caption{Illustration of the construction preceding Lemma~\ref{lem:sign-preserving2}:
If $r_i^\prime$ and $b_i^\prime$ are of different colour, then $e_i^r$ and $e_i^b$
must necessarily both belong to every \som\ ${\red\mu^\prime}\cup{\blue\nu^\prime}$
and thus can be removed (after accounting for their weights, of course) together with
the vertices $r_i^\prime$ and $b_i^\prime$. The resulting situation ``locally'' looks
like the original graph $G$, with coloured vertices $r_i$ (of the same colour as $b_i^\prime$)
and $b_i$  (of the same colour as $r_i^\prime$).
}
\label{fig:edge-condensation2}
\input graphics/edge-condensation2
\end{figure}

If $r_i^\prime$ and $b_i^\prime$ are of the \EM{same} colour (say red), then for any superposition of matchings ${\red\mu}^\prime\cup{\blue\nu}^\prime$
in $B^\prime$ 
\bit
\item
the blue subgraph $\subgraph{G^\prime}{\thereds^\prime}$ (which contains the 
\EM{blue} matching ${\blue\nu}^\prime$)
%  clearly is a matching of
%$\subgraph{G^\prime}{\setof{r_i^\prime,b_i^\prime}}$, which
``locally looks like''
the original graph $G$ with edge $e_i$ removed (since this is ``locally equivalent''
to removing the two red
vertices $r_i^\prime$ and $b_i^\prime$, see the lower left picture in
\figref{fig:edge-condensation3}),
\item
and in the red subgraph  $\subgraph{G^\prime}{\theblues^\prime}$ (which contains the 
red matching ${\red\mu}^\prime$)
we may re--replace the path of length 3 $\pas{r_i,r_i^\prime,b_i^\prime,b_i}$
in $G^\prime$ by the ``original'' edge $e_i=\setof{r_i,b_i}$ in $G$: This operation changes
the corresponding Pfaffian only by a sign--factor according to  the considerations preceding
\eqref{eq:ciucu1996}, and the resulting red subgraph ``locally looks like'' the original graph $G$
(see the lower right picture in \figref{fig:edge-condensation3}).
\eit
% \eit
\begin{figure}
\caption{Illustration of the construction preceding Lemma~\ref{lem:sign-preserving2}:
If $r_i^\prime$ and $b_i^\prime$ are of the same colour (say red), then the construction
shown in \figref{fig:edge-condensation} can be reversed in the ``red subgraph'' (see the
right lower picture). This operation introduces a sign factor according to the considerations preceding
\eqref{eq:ciucu1996}.
In the ``blue subgraph'' the vertices $r_i^\prime$ and $b_i^\prime$ simply are missing
(see the left lower picture).
}
\label{fig:edge-condensation3}
\input graphics/edge-condensation3
\end{figure}

Applying these simple considerations to \EM{all}
pairs $\pas{r_i^\prime,b_i^\prime}$, $i=1,\dots, k$ in $G^\prime$, we see that superpositions of matchings in
$B^\prime=\bicoloured{G^\prime}{\theblues^\prime}{\thereds^\prime}$
% any bicoloured graph $B$ derived
%from $G^\prime$
%(with coloured vertices $\thecoloureds=\setof{r_1^\prime,b_1^\prime,\dots, r_k^\prime,b_k^\prime}$)
are in bijection with superpositions of matchings in a certain bicoloured graph derived
from $G$, which we now describe:
%. So if we properly account for sign--changes and weight--changes, we obtain
%another Pfaffian identity:
% in the sense
%that all the ``additional vertices''
%$r_i^\prime$ and $b_i^\prime$ can be removed without modifying the corresponding Pfaffians.
%To be more precise:
%, consider the
%following partition of the set $\setof{1,\dots,k}$ in disjoint subsets, which is
%induced by the colouring of $\thecoloureds=\thereds^\prime\disjuni\theblues^\prime$:
%\bit
%\item $\bb\defeq\setof{i: r_i^\prime\text{ and } b_i^\prime\text{ are both blue}}$,
%\item $\rr\defeq\setof{i: r_i^\prime\text{ and } b_i^\prime\text{ are both red}}$,
%% \item $\rb\defeq\setof{i: r_i^\prime\text{ and } b_i^\prime\text{ are of different colour}}$.
%\item $\rb\defeq\setof{i: r_i^\prime\text{ is red, } b_i^\prime\text{ is blue}}$,
%\item $\br\defeq\setof{i: r_i^\prime\text{ is blue, } b_i^\prime\text{ is red}}$.
%\eit

Define the sets of blue/red \EM{vertices}
in $G$, whose ``partner with the same subscript'' has a
different colour (see \figref{fig:edge-condensation2}):
%$\theblues^{\prime\prime}$ (new blue vertices)
%and $\thereds^{\prime\prime}$ (new red vertices) of \EM{vertices}
%in $G$ by
\begin{align}
\theblues^{\prime\prime} &\defeq
\setof{r_i:\; r_i^\prime\in\thereds^\prime\and b_i^\prime\in\theblues^\prime}
\cup
\setof{b_i:\; b_i^\prime\in\thereds^\prime\and r_i^\prime\in\theblues^\prime}
%\cup
%\setof{b_i, r_i:\; b_i^\prime\in\theblues^\prime\and r_i^\prime\in\theblues^\prime}
, \notag \\
\thereds^{\prime\prime} &\defeq
\setof{b_i:\; r_i^\prime\in\thereds^\prime\and b_i^\prime\in\theblues^\prime}
\cup
\setof{r_i:\; b_i^\prime\in\thereds^\prime\and r_i^\prime\in\theblues^\prime}.
\label{eq:primeprime}
%\cup
%\setof{b_i, r_i:\; b_i^\prime\in\thereds^\prime\and r_i^\prime\in\thereds^\prime}.
%, \\
%\theblues^{\prime\prime} &\defeq
%\setof{r_i:\; r_i^\prime\text{ is red, but } b_i^\prime\text{ is blue}}
%\cup
%\setof{b_i:\; b_i^\prime\text{ is red, but } r_i^\prime\text{ is blue}}
%, \\
%\thereds^{\prime\prime} &\defeq
%\setof{r_i:\; r_i^\prime\text{ is blue, but } b_i^\prime\text{ is red}}
%\cup
%\setof{b_i:\; b_i^\prime\text{ is blue, but } r_i^\prime\text{ is red}}.
\end{align}
Define the sets of blue/red \EM{edges} in $G$, where both end--vertices are of the
same colour (see \figref{fig:edge-condensation3}):
% $\theblues_e$ and $\thereds_e$ of \EM{edges} in $G$ by
\begin{align}
\thereds_e &\defeq
\setof{e_i:\; b_i^\prime\in\thereds^\prime\and r_i^\prime\in\thereds^\prime}, \notag \\
\theblues_e &\defeq
\setof{e_i:\; b_i^\prime\in\theblues^\prime\and r_i^\prime\in\theblues^\prime}.\label{eq:edges}
%, \\
%\theblues_e &\defeq
%\setof{e_i:\; r_i^\prime\text{ and } b_i^\prime\text{ are both blue}}, \\
%\thereds_e &\defeq
%\setof{e_i:\; r_i^\prime\text{ and } b_i^\prime\text{ are both red}}.
\end{align}
%(Note that these lenghty definitions simply follow from the construction
%illustrated by Figures~\ref{fig:edge-condensation2} and \ref{fig:edge-condensation3}.)

Finally, define the set of vertices incident with edges from $\thereds_e\cup\theblues_e$:
\begin{equation*}
Z\defeq\setof{b_i, r_i:\;b_i \text{ and } r_i\text{ are of the same colour}}
=\pas{\bigcup_{e\in\thereds_e}e}\cup\pas{\bigcup_{e\in\theblues_e}e}.
\end{equation*}

For an arbitrary graph $H$ and subsets $V^\prime\subseteq\vertices\of{H}$ and
$E^\prime\subseteq\edges\of{H}$
introduce the notation
$$\subgraphedges{H}{V^\prime}{E^\prime}$$
 for the
graph $\subgraph{H}{V^\prime}$ with all % from which all
edges in $E^\prime$ removed.
Then the above reasoning amounts to the following Pfaffian identity:
\begin{multline}
\label{eq:lemma-yyz}
\pfaffian\of{\subgraph{G^\prime}{\theblues^\prime}}\cdot
\pfaffian\of{\subgraph{G^\prime}{\thereds^\prime}}
= \\ \pas{-1}^{\cardof{Z}+\SETSUM{Z}{V}}
\pas{\prod_{e_i\not\in\pas{\theblues_e\cup\thereds_e}}\weight\of{e_i}}\cdot
\pfaffian\of{\subgraphedges{G}{\theblues^{\prime\prime}}{\theblues_e}}\cdot
\pfaffian\of{\subgraphedges{G}{\thereds^{\prime\prime}}{\thereds_e}}.
\end{multline}
%Now apply Lemma~\ref{lem:sign-preserving} to $G^\prime$, and translate the
%identity \eqref{eq:sign-preserving} according to \eqref{eq:lemma-yyz}.
These considerations lead to 
% immediately proves
the following assertion:
\begin{lem}
\label{lem:sign-preserving2}
Assume the situation of Lemma~\ref{lem:sign-preserving} with the additional
condition that there is \EM{no} vertex between $r_i$ and $b_i$ in the ordered
set of vertices $\vertices\of G$ ($G=K_V$). For an arbitrary subset
$I\subseteq\setof{1,\dots,k}$
denote by $\complement{I}$ the subset
$\setof{1,\dots,k}\setminus I$ and 
introduce the following ``template notation'':
% \begin{align*}
%\complement{I}  &\defeq \setof{1,\dots,k}\setminus I, \\
%\redindexed{I}  &\defeq \setof{r_i:\;i\in I}, \\
%\blueindexed{I} &\defeq \setof{b_i:\;i\in I} , \\
%\edgeindexed{I} &\defeq \setof{e_i:\;i\in I} .
$$ \indexed{x}{I} \defeq \setof{x_i:\;i\in I}, $$
%\end{align*}
where $x$ may be any symbol from the set $\setof{r,r^\prime,b,b^\prime, e,e^\prime}$.

%Then we have for all (fixed) subsets $X\subseteq\thereds$:
%\begin{multline}
%% \label{eq:sign-preserving2}
%\sum_{W\subseteq\theblues}
%	\pfaffian\of{\pas{\thereds\cup\thewhites}\cup W}
%	\cdot
%	\pfaffian\of{\pas{\theblues\cup\thewhites}\setminus W} \\
%=
%\sum_{V\subseteq\theblues}
%	\pfaffian\of{\pas{\pas{\thereds\cup\thewhites}\cup V} \setminus X}
%	\cdot
%	\pfaffian\of{\pas{\pas{\theblues\cup\thewhites}\setminus V}\cup X}.
%\end{multline}

Then we have for all \EM{fixed} subsets $B\subseteq\myrange{k}$:
\begin{multline}
\label{eq:sign-preserving2}
\sum_{R\subseteq\myrange{k}}
\pas{\prod_{e_i\in\indexed{e}{\complement{R}}}\weight\of{e_i}}
\cdot
\pfaffian\of{\subgraph{G}{\indexed{r}{\complement{R}}}}
\cdot
\pfaffian\of{\subgraphedges{G}{
\indexed{b}{\complement{R}}
}{\indexed{e}{R}}}
=\\
\sum_{R\subseteq\myrange{k}}
\pas{\prod_{e_i\in\pas{\indexed{e}{\complement{\symdef{R}{B}}}}}\weight\of{e_i}}
\times \\
\pfaffian\of{\subgraphedges{G}{%
\pas{\indexed{r}{\complement{B}\cap\complement{R}} \cup\indexed{b}{R\cap B}}
}{\indexed{e}{B\setminus R}}}
\cdot
\pfaffian\of{\subgraphedges{G}{
\pas{\indexed{b}{\complement{R}\cap\complement{B}} \cup \indexed{r}{B\cap R}}
}{\indexed{e}{R\setminus B}}}
\end{multline}
This identity is of the \wpit.
\end{lem}
\begin{proof}
Consider the graph $G^\prime$ as defined in the considerations preceding
Lemma~\ref{lem:sign-preserving2}.
It is obvious that every bicoloured path in $G$ connecting two vertices $b_i$ and $b_j$
(or $r_i$ and $r_j$) corresponds bijectively to a bicoloured path in $G^\prime$
% (of the same weight??????????)
 connecting the
vertices $b_i^\prime$ and $b_j^\prime$ (or  $r_i^\prime$ and $r_j^\prime$), which
thus also has weight zero. Hence the
weight function of $G^\prime$ is also planar, and we may apply
Lemma~\ref{lem:sign-preserving} to $G^\prime$.
%, and all we have to is to
%``translate''
%\eqref{eq:sign-preserving} properly.

Let $R,B\subseteq\setof{1,\dots,k}$ such that % $V\subseteq\therights^\prime$ and $X\subseteq\thelefts^\prime$
% be given by index sets $S,C\subseteq\setof{1,\dots,k}$, i.e.,
$V=\indexed{b^\prime}{R}\subseteq\therights^\prime$ and
as $X=\indexed{r^\prime}{B}\subseteq\thelefts^\prime$. Then we may write
% also have  the translation
\begin{equation}
\label{eq:translate1}
\pfaffian\of{\pas{\pas{\thelefts^\prime\cup\thewhites}\cup V} \setminus X}
	\cdot
	\pfaffian\of{\pas{\pas{\therights^\prime\cup\thewhites}\setminus V}\cup X}
=
\pfaffian\of{\subgraph{G^\prime}{\theblues^\prime}}
\pfaffian\of{\subgraph{G^\prime}{\thereds^\prime}},
\end{equation}
where
$\theblues^\prime = \indexed{r}{B} \cup \indexed{b}{\complement{R}}$ and
$\thereds^\prime  = \indexed{b}{R} \cup \indexed{r}{\complement{B}}$.
%Consider two arbitrary subsets $R, B\subseteq\setof{1\dots,k}$, which define a partition
%of the set of coloured vertices $\setof{r_1^\prime,b_1^\prime,\dots,r_k^\prime,b_k^\prime}$ of $G^\prime$ into blue and red
%subsets
%\begin{align*}
%\theblues^\prime &= \indexed{r}{B} \cup \indexed{b}{\complement{R}}, \\
%\thereds^\prime  &= \indexed{b}{R} \cup \indexed{r}{\complement{B}}.
%\end{align*}
%(These definitions can easily be interpreted as follows:
%Imagine that \EM{all} vertices in $\therights^\prime$ are ``initially'' blue, and
%$\indexed{b^\prime}{R}\subseteq\therights^\prime$ is a subset of vertices whose initial colour is
%changed to red, and that \EM{all} vertices in $\thelefts^\prime$ are ``initially'' red, and
%$\indexed{r^\prime}{B}\subseteq\thelefts^\prime$ is a subset of vertices whose initial colour is
%changed to blue.)

The definitions of $\theblues^{\prime\prime}$ and $\thereds^{\prime\prime}$
(according to \eqref{eq:primeprime}) and of $\theblues_e$ and $\thereds_e$
(according to \eqref{eq:edges}) read:
\begin{align*}
\theblues^{\prime\prime} &=
	\indexed{r}{\complement{B}\cap\complement{R}} \cup
	\indexed{b}{R\cap B}, \\
\thereds^{\prime\prime}  &= 
	\indexed{b}{\complement{R}\cap\complement{B}} \cup
	\indexed{r}{B\cap R}, \\
%\end{align*}
%And the definitions of $\theblues_e$ and $\thereds_e$
%(according to \eqref{eq:edges}) then read:
%\begin{align*}
\theblues_e &=
	\indexed{e}{B\setminus R}, \\
\thereds_e  &= 
	\indexed{e}{R\setminus B}.
\end{align*}

Then \eqref{eq:lemma-yyz} translates to:
\begin{multline}
\label{eq:translate2}
\pfaffian\of{\subgraph{G^\prime}{\theblues^\prime}}
\pfaffian\of{\subgraph{G^\prime}{\thereds^\prime}}
= \pas{\prod_{e_i\in\pas{\indexed{e}{\complement{\symdef{R}{B}}}}}\weight\of{e_i}}
\times \\
\pfaffian\of{\subgraphedges{G}{%
\pas{\indexed{r}{\complement{B}\cap\complement{R}} \cup\indexed{b}{R\cap B}}
}{\indexed{e}{B\setminus R}}}\cdot
\pfaffian\of{\subgraphedges{G}{
\pas{\indexed{b}{\complement{R}\cap\complement{B}} \cup \indexed{r}{B\cap R}}
}{\indexed{e}{R\setminus B}}}.
\end{multline}
(Note that there is no sign--change here, since there are no vertices between
$r_i$ and $b_i$ by assumption.)

Putting together \eqref{eq:translate1} and \eqref{eq:translate2} gives the translation from
\eqref{eq:sign-preserving} to \eqref{eq:sign-preserving2}.
%where 
%$\theblues^\prime = \indexed{r}{C} \cup \indexed{b}{\complement{S}}$
%and
%$\thereds^\prime  = \indexed{b}{S} \cup \indexed{r}{\complement{C}}$.
% \end{align*}
\end{proof}

This yields immediately the following generalization of a result by Yan, Yeh and Zhang
\cite[Theorem~3.2]{yan-yeh-zhang:2005} (here, we write the summands for $\complement{R}$ instead of $R$ when summing over all $R\subseteq\myrange{k}$, to make the
equivalence with \cite[equation (8)]{yan-yeh-zhang:2005} more transparent):
\begin{cor}
\label{cor:yyz-theorem}
Let $G$ be a \EM{planar} graph with $k$ independent edges $e_i=\setof{a_i,e_i}$,
$i=1,\dots, k$, in the
boundary of some face $f$  of $G$, such that the vertices $a_1,b_1,\dots,a_k,b_k$ 
appear in that cyclic order in $f$.

Then we have for all \EM{fixed} subsets $B\subseteq\myrange{k}$:
\begin{multline*}
\sum_{{R}\subseteq \myrange{K}}
\pas{\prod_{e_i\in\indexed{e}{{R}}}\weight\of{e_i}}
\cdot
\nofmatchings{\subgraph{G}{\indexed{r}{{R}}}}
\cdot
\nofmatchings{\subgraphedges{G}{
\indexed{b}{{R}}
}{\indexed{e}{\complement{R}}}}
= \\
\sum_{{R}\subseteq\myrange{k}}
\pas{\prod_{e_i\in\pas{\indexed{e}{{\symdef{R}{B}}}}}\weight\of{e_i}}
\cdot
\nofmatchings{\subgraphedges{G}{%
\pas{\indexed{r}{\complement{B}\cap{R}} \cup\indexed{b}{\complement{R}\cap B}}
}{\indexed{e}{B\setminus\complement{R}}}}
\cdot
\nofmatchings{\subgraphedges{G}{
\pas{\indexed{b}{{R}\cap\complement{B}} \cup \indexed{r}{B\cap\complement{R}}}
}{\indexed{e}{\complement{R}\setminus B}}}
\end{multline*}

%Moreover, assume
%that Let $\thereds=\setof{a_1,\dots,a_k}$ and $\theblues=\setof{b_1,\dots,b_k}$. Then we have for every fixed subset $X\subseteq\thereds$
%\begin{equation*}
%\sum_{W\subseteq \theblues}
%	\nofmatchings{\subgraph{G}{\pas{\theblues\setminus W}}}
%	\nofmatchings{\subgraph{G}{\pas{\thereds \cup {W}}}}
%=
%\sum_{V\subseteq \theblues}
%	\nofmatchings{\subgraph{G}{\pas{\pas{\theblues\setminus V} \cup X}}}
%	\nofmatchings{\subgraph{G}{\pas{\pas{\thereds\cup V} \setminus X}}}.
%\end{equation*}
\end{cor}
\begin{proof}
The statement follows from \eqref{eq:sign-preserving2} by the Kasteleyn--Percus
method.
\end{proof}

\bibliography{paper}

\end{document}